\documentclass[11pt, reqno]{amsart}
\usepackage{mathrsfs}
\usepackage[active]{srcltx}
\usepackage{mathrsfs}
\usepackage{array}
\usepackage{multirow}
\usepackage{float}
\usepackage{caption}
\usepackage{yfonts}
\DeclareMathAlphabet{\mathpzc}{OT1}{pzc}{m}{rm}

\usepackage{color}
\usepackage[unicode]{hyperref}
\hypersetup{
	colorlinks = true,%
	citecolor = [rgb]{0.0,0.0,0.9},
	filecolor=black,%
	linkcolor = [rgb]{0.65,0.0,0.0},%
	anchorcolor = red,
	pagecolor = red,
	urlcolor= [rgb]{0.65,0.0,0.0}
}

\newtheorem{rem}{Remark}[section]

\numberwithin{equation}{section}

\usepackage{graphicx}

\usepackage[left=2.85cm, right=2.85cm, top=2.0cm, bottom=2.00cm]{geometry}

\newtheorem{proposition}{Proposition}[section]
\newtheorem{theorem}{Theorem}[section]

\newtheorem{remark}{Remark}[section]

\newcommand{\ds}{\displaystyle}

\begin{document}
\title[Clamped plates in nonpositively curved spaces]
{
	Fundamental tones of clamped plates in nonpositively curved spaces}

\vspace{-0.9cm}
\author{Alexandru Krist\'aly}
\address{Institute of Applied Mathematics, \'Obuda
	University, 
	Budapest, Hungary \& Department of Economics, Babe\c s-Bolyai University, Cluj-Napoca, Romania  
}
\email{alex.kristaly@econ.ubbcluj.ro; kristaly.alexandru@nik.uni-obuda.hu
}

\thanks{Research supported by the National Research, Development and Innovation Fund of Hungary, financed under the K$\_$18 funding scheme, Project No.  127926.}

\keywords{Lord Rayleigh's isoperimetric problem; fundamental tone; clamped plate; nonpositive curvature;   hypergeometric function.}

\subjclass[2000]{Primary: 35P15, 53C21, 35J35, 35J40.
}

\begin{abstract} 
	We study Lord Rayleigh's problem for clamped plates on an arbitrary $n$-dimen\-sional $(n\geq 2)$ Cartan-Hada\-mard manifold $(M,g)$  with sectional curvature $\textbf{K}\leq -\kappa^2$ for some $\kappa\geq 0.$ We first prove a McKean-type spectral gap estimate, i.e.  the fundamental tone of any domain in $(M,g)$ is universally bounded  from below by $\frac{(n-1)^4}{16}\kappa^4$ whenever the $\kappa$-Cartan-Hada\-mard conjecture holds on $(M,g)$, e.g. in 2-\ and 3-dimensions due to Bol (1941) and Kleiner (1992), respectively. In  2- and 3-dimensions we prove  sharp isoperimetric inequalities for sufficiently small clamped plates, i.e.  the fundamental tone of any domain in $(M,g)$ of volume $v>0$ is not less than the corresponding fundamental tone of a geodesic ball of the same volume $v$ in the space of constant curvature $-\kappa^2$ provided that  $v\leq c_n/\kappa^n$ with $c_2\approx 21.031$ and  $c_3\approx 1.721$, respectively. In particular, Rayleigh's  problem in  Euclidean spaces  resolved by  Nadirashvili (1992) and Ashbaugh and Benguria (1995) appears as a limiting case in our setting (i.e.  $\textbf{K}\equiv\kappa=0$). Sharp asymptotic estimates of the fundamental tone of small and large geodesic balls of low-dimensional hyperbolic spaces are also given. The sharpness of our results requires the validity of the $\kappa$-Cartan-Hada\-mard conjecture (i.e.  sharp isoperimetric inequality on $(M,g)$) and peculiar properties of the Gaussian hypergeometric function, both valid only in dimensions 2 and 3; nevertheless, some nonoptimal estimates of the fundamental tone of arbitrary clamped plates are also provided in high-dimensions.  As an application, by using the sharp isoperimetric inequality for small clamped hyperbolic discs, we give necessarily and sufficient conditions for the existence of a nontrivial solution to an elliptic PDE involving the biharmonic Laplace-Beltrami operator.  
\end{abstract}

\maketitle
\vspace{-0.7cm}
\begin{quotation}
	Dedicated to Professor Biagio Ricceri on the occasion of his 65th anniversary.
\end{quotation}

\vspace{-0.3cm}
\section{Introduction and Main Results}


Let $\Omega\subset \mathbb R^n$ be a bounded domain $(n\geq 2)$, and consider the eigenvalue problem 
\begin{equation}\label{R}
\left\{ \begin{array}{lll}
\Delta^2 u=\Gamma u &\mbox{in} &  \Omega, \\
u=|\nabla u|=0 &\mbox{on} &  \partial \Omega,
\end{array}\right.
\end{equation}
 associated with the vibration of a clamped plate. The lowest/principal eigenvalue for (\ref{R}) -- the fundamental tone of the clamped plate -- can be characterized in a variational way by 
   \begin{equation}\label{variational-charact-0}
   \Gamma_0(\Omega)=\inf_{u\in W_0^{2,2}(\Omega)\setminus \{0\}}\frac{\displaystyle \int_{\Omega}(\Delta u)^2 {\rm d}x}{\displaystyle \int_{\Omega}u^2 {\rm d}x}.
   \end{equation} 
   The minimizer of (\ref{variational-charact-0}) in the plane describes the vibration of a homogeneous thin plate $\Omega\subset \mathbb R^2$ whose boundary is clamped, while the frequency of vibration of the plate $\Omega$ is proportional to  $\Gamma_0(\Omega)^\frac{1}{2}.$
The famous conjecture of Lord Rayleigh \cite[p.382]{Rayleigh} -- formulated initially for planar domains in 1894 -- states that 
\begin{equation}\label{Rayleigh-eredeti}
\Gamma_0(\Omega)\geq \Gamma_0(\Omega^\star)=\mathfrak  h^4_{\nu}\left(\frac{\omega_n}{|\Omega|}\right)^{4/n},
\end{equation}
where $\Omega^\star\subset \mathbb R^n$ is a ball with the same measure as $\Omega$, with equality if and only if $\Omega$ is  a ball. Hereafter, $\nu=\frac{n}{2}-1$, $\omega_n=\pi^{n/2}/\Gamma(1+n/2)$  is the volume of the unit Euclidean ball, while $\mathfrak  h_\nu$ is the 
 first positive critical point of $\frac{J_\nu}{I_\nu}$, where $J_\nu$ and $I_\nu$ stand for the  Bessel and modified Bessel functions of first kind, respectively. 
 
 Assuming that the eigenfunction corresponding to $\Gamma_0(\Omega)$ is sign-preserving over a simply connected domain $\Omega\subset \mathbb R^2$,  Szeg\H o \cite{Szego} proved (\ref{Rayleigh-eredeti}) in the early fifties.\ As one can deduce from his paper's text, his belief on the constant-sign first eigenfunction corresponding to $\Gamma_0(\Omega)$ has been based on the second-order membrane problem (called as the Faber-Krahn problem). It turned out shortly that his expectation perishes due to the construction of Duffin \cite{Duffin} on strip-like domains and Coffman, Duffin and Shaffer \cite{CDS} on ring-shaped clamped plate, localizing \textit{nodal lines} of vibrating plates. While the membrane problem involves only the Laplacian, the clamped plate problem requires the presence of the fourth order bilaplacian operator; as we know nowadays, fourth order equations are lacking general maximum/comparison principles which is unrevealed  in Szeg\H o's pioneering approach. In fact, stimulated by the papers  \cite{Duffin} and \cite{CDS}, several scenarios are described for nodal domains of clamped plates, see e.g. Bauer and Reiss \cite{Bauer}, Coffman \cite{Coffman}, Grunau and Robert \cite{Grunau},  from which the main edification is that eigenfunctions corresponding to (\ref{variational-charact-0})  may change their sign. 
 

  
  In order to handle the presence of possible nodal domains, Talenti \cite{Talenti} developed a Schwarz-type rearrangement method  on domains where the first eigenfunction corresponding to (\ref{variational-charact-0}) has both positive and negative parts. In this way, a decomposition of  (\ref{variational-charact-0}) into a two-ball minimization problem arises which provided a  nonoptimal estimate in (\ref{Rayleigh-eredeti}); in fact, instead of (\ref{Rayleigh-eredeti}), Talenti proved that $\Gamma_0(\Omega)\geq d_n\Gamma_0(\Omega^\star)$ where the dimension-depending constant $d_n$ has the properties $\frac{1}{2}\leq d_n<1$ for every $n\geq 2$ and $\ds\lim_{n\to \infty}d_n=\frac{1}{2}.$
  
  By a careful improvement of Talenti's two-ball minimization argument,  Rayleigh's conjecture has been proved in its full generality for $n=2$ by Nadirashvili \cite{Nad-0, Nad}. Further modifications of some arguments from  the papers \cite{Nad} and \cite{Talenti} allowed to Ashbaugh and Benguria \cite{A-B} to prove Rayleigh's conjecture for $n=3$ (and $n=2$) by exploring fine properties of Bessel functions. Roughly speaking, for $n\in \{2,3\}$, the two-ball minimization problem reduces to only one ball (the other ball disappearing), while in higher dimensions the 'optimal' situation appears for two identical balls which provides a nonoptimal estimate for $\Gamma_0(\Omega)$.  Although asymptotically sharp estimates are provided  by Ashbaugh and Laugesen \cite{A-L} for $\Gamma_0(\Omega)$ in high-dimensions, i.e.  $\ds\Gamma_0(\Omega)\geq D_n\Gamma_0(\Omega^\star)$ where $0.89<D_n<1$ for every $n\geq 4$ with $\ds\lim_{n\to \infty}D_n=1$, the conjecture is still open for $n\geq 4$. Very recently, Chasman and Langford \cite{CL1,CL2} provided certain Ashbaugh-Laugesen-type results 
  in Euclidean spaces endowed with a log-convex/Gaussian density, by proving that $\Gamma_w(\Omega)\geq \tilde C\Gamma_w(\Omega^\star)$, where the constant $\tilde C\in (0,1)$ depends  on  the  volume of $\Omega$ and  dimension $n\geq 2$, while  $\Gamma_w(\Omega)$ and $\Gamma_w(\Omega^\star)$ denote the fundamental tones of the clamped plate with respect to the corresponding density function $w$.
   
 Interest in the clamped plate problem on curved spaces was also increased 
in recent years. One of the most central problems is to establish  Payne-P\'olya-Weinberger-Yang-type inequalities for the eigenvalues of the   problem   
   \begin{equation}\label{CP-problem}
   \left\{ \begin{array}{lll}
   \Delta_g^2 u=\Gamma u &\mbox{in} &  \Omega, \\
   u=\frac{\partial u}{\partial \textbf{n}}=0 &\mbox{on} &  \partial \Omega,
   \end{array}\right.
   \end{equation}
   where $\Omega$ is a bounded domain  in an $n$-dimensional Riemannian manifold $(M,g)$,    $\Delta_g^2$ stands for the biharmonic Laplace-Beltrami operator on $(M,g)$ and $\frac{\partial }{\partial \textbf{n}}$ is the outward normal   derivative on $\partial \Omega$, respectively; see e.g. Chen, Zheng and Lu \cite{Chen-Zheng-Lu}, Cheng,  Ichikawa and Mametsuka \cite{Cheng-Ichikawa-Mametsuka}, Cheng and Yang \cite{Cheng-Yang-0, Cheng-Yang-1, Cheng-Yang-2}, Wang and Xia \cite{Wang-Xia}. 
     Instead of (\ref{variational-charact-0}), one  naturally considers the \textit{fundamental tone of} $\Omega\subset M$ by  
   \begin{equation}\label{variational-charact}
   \Gamma_g(\Omega):=\Gamma_{g,n}(\Omega)=\inf_{u\in W_0^{2,2}(\Omega)\setminus \{0\}}\frac{\displaystyle \int_{\Omega}(\Delta_g u)^2 {\rm d}v_g}{\displaystyle \int_{\Omega}u^2 {\rm d}v_g},
   \end{equation} 
   where ${\rm d}v_g$ denotes the canonical measure on $(M,g)$,  and $W_0^{2,2}(\Omega)$ is the usual Sobolev space on $(M,g)$, see Hebey \cite{Hebey}; in fact, it turns out that  $\Gamma_g(\Omega)$ is the first eigenvalue to (\ref{CP-problem}).   Due to the Bochner-Lichnerowicz-Weitzenb\"ock formula, the Sobolev space  $H_0^{2}(\Omega)=W_0^{2,2}(\Omega)$ is a proper choice for (\ref{CP-problem}),  see Proposition \ref{prop-1} for details.  
 
 To the best of our knowledge, no results -- comparable to (\ref{Rayleigh-eredeti}) -- are available in the literature concerning Lord Rayleigh's problem for clamped plates on curved structures. Accordingly, the main purpose of the present paper is to identify those geometric and analytic properties which reside in Lord Rayleigh's problem for clamped plates on \textit{nonpositively curved spaces}. To develop our results, the  geometric context is provided by an $n$-dimensional $(n\geq 2)$ Cartan-Hadamard manifold $(M,g)$ (i.e.  simply connected, complete Riemannian manifold with nonpositive sectional curvature). Having this framework, we recall McKean's spectral gap estimate for membranes which is closely related to (\ref{variational-charact}); namely, in an $n$-dimensional Cartan-Hadamard manifold $(M,g)$  with sectional curvature ${{\rm \bf K}}\leq -\kappa^2$ for some $\kappa> 0,$  the 
 principal frequency of any membrane $\Omega\subset M$ can be estimated as
 \begin{equation}\label{variational-charact-1}
 \gamma_g(\Omega):=\inf_{u\in W_0^{1,2}(\Omega)\setminus \{0\}}\frac{\displaystyle \int_{\Omega}|\nabla_g u|^2 {\rm d}v_g}{\displaystyle \int_{\Omega}u^2 {\rm d}v_g}\geq \frac{(n-1)^2}{4}\kappa^2;
 \end{equation} 
 in addition, (\ref{variational-charact-1}) is sharp on the $n$-dimensional hyperbolic space $(\mathbb H_{-\kappa^2}^n,g_\kappa)$  of constant curvature $-\kappa^2$ in the sense that $\gamma_{g_\kappa}(\Omega)\to \frac{(n-1)^2}{4}\kappa^2$ whenever $\Omega$ tends to $\mathbb H_{-\kappa^2}^n$, see McKean \cite{McKean}. 
 
 Before to state our results, we fix some notations. If $\kappa\geq 0$, let $N^n_{\kappa}$ be the $n$-dimensional space-form  with constant sectional curvature $-\kappa^2$, i.e.  $N^n_{\kappa}$ is either the hyperbolic space $\mathbb H_{-\kappa^2}^n$ when $k>0$, or the Euclidean space $\mathbb R^n$ when $\kappa=0.$ Let $B_\kappa(L)$ be the geodesic ball of radius $L>0$ in $N^n_{\kappa}$ and if $\Omega\subset N^n_{\kappa}$, we denote by $\Gamma_\kappa(\Omega)$ the corresponding value from (\ref{variational-charact}). 
 By convention, we consider $1/0=+\infty$ and as usual, $V_g(S)$ denotes the Riemannian volume of $S\subset M$.

 Our first result provides a fourth order counterpart of McKean's spectral gap estimate, which requires the validity of the \textit{$\kappa$-Cartan-Hadamard conjecture on} $(M,g)$; the latter is nothing but the sharp isoperimetric inequality on $(M,g)$, which is valid e.g.\ on hyperbolic spaces of any dimension as well as on generic 2- and 3-dimensional Cartan-Hadamard manifolds with sectional curvature ${{\rm \bf K}}\leq -\kappa^2$ for some $\kappa\geq 0,$ see \S \ref{CH-conjecture}.

 \begin{theorem}\label{main-theorem}
 	Let  $(M,g)$ be an $n$-dimensional Cartan-Hadamard manifold with sectional curvature ${{\rm \bf K}}\leq -\kappa^2$ for some $\kappa\geq 0,$ which verifies the $\kappa$-Cartan-Hadamard conjecture.  If  $\Omega\subset M$  is a bounded domain with smooth boundary then 
 	\begin{equation}\label{lambda-becsles}
 	\Gamma_g(\Omega)\geq  \frac{(n-1)^4}{16}\kappa^4.
 	\end{equation}
 	Moreover, for $n\in \{2,3\}$, relation {\rm (\ref{lambda-becsles})} is sharp in the sense that 
 \begin{equation}\label{hatar-vegtelen}
 \Gamma_\kappa(N^n_{\kappa}):=\lim_{L\to \infty}\Gamma_\kappa(B_\kappa(L))= \frac{(n-1)^4}{16}\kappa^4.
 \end{equation}

 \end{theorem}

 Clearly, Theorem \ref{main-theorem} is relevant only for $\kappa>0$ (as  (\ref{lambda-becsles}) and (\ref{hatar-vegtelen}) trivially hold for $\kappa=0$). Moreover, if  $n\in \{2,3\}$ and $\kappa>0$, and  $\Gamma_\kappa^l(\Omega)$ denotes the $l^{\rm th}$ eigenvalue of {\rm (\ref{CP-problem})} on $\Omega\subset \mathbb H_{-\kappa^2}^n$, then making use of  (\ref{hatar-vegtelen})  and a Payne-P\'olya-Weinberger-Yang-type universal inequality on $\mathbb H_{-\kappa^2}^n$,
  it turns out that   \begin{equation}\label{magasabbbakkk}
  \Gamma_\kappa^l(\mathbb H_{-\kappa^2}^n):=\lim_{L\to \infty}\Gamma_\kappa^l(B_\kappa(L))=\frac{(n-1)^4}{16}\kappa^4\ \ {\rm for\ all}\ \ l\in \mathbb N. 
  \end{equation}
 In particular, (\ref{magasabbbakkk}) confirms a claim of Cheng and Yang \cite[Theorem 1.4]{Cheng-Yang-1} for $n\in \{2,3\}$, where the authors assumed (\ref{hatar-vegtelen}) itself in order to derive (\ref{magasabbbakkk}). In fact, one can expect the validity of (\ref{magasabbbakkk}) for any $n\geq 2$ but some technical difficulties prevent the proof in high-dimensions;  for details, see \S \ref{hatar-helyyzettt}.

 Actually, Theorem \ref{main-theorem} is just a byproduct of a general argument needed to prove the main result of our paper (for its statement, we recall that  $\mathfrak  h_\nu$ is the 
 first positive critical point of $\frac{J_\nu}{I_\nu}$ and $\nu=\frac{n}{2}-1$):

 \begin{theorem}\label{masodik-fotetel}  	Let $n\in \{2,3\}$ and  $(M,g)$ be an $n$-dimensional Cartan-Hadamard manifold with sectional curvature ${{\rm \bf K}}\leq -\kappa^2$ for some $\kappa\geq 0,$  $\Omega\subset M$  be a bounded domain with smooth boundary and volume  $V_g(\Omega)\leq \frac{c_n}{\kappa^n}$ with  $c_2\approx 21.031$ and $c_3\approx 1.721.$
 	If $\Omega^\star\subset N^n_\kappa$ is a geodesic ball 
 	verifying $V_g(\Omega)=V_\kappa(\Omega^\star)$	then 
 	\begin{equation}\label{A_B_egyenlotlenseg}
 	\Gamma_g(\Omega)\geq \Gamma_\kappa(\Omega^\star),
 	\end{equation}
 	with  equality in {\rm (\ref{A_B_egyenlotlenseg})} if and only if $\Omega$ is isometric to $\Omega^\star.$ Moreover, 
 	\begin{equation}\label{asymptotics}
 	\Gamma_\kappa(B_\kappa(L))\sim\left({\frac{(n-1)^2}{4}{\kappa^2}+\frac{\mathfrak  h_\nu^2}{L^2}}\right)^2\ \ \mbox{as}\ \ L\to 0.
 	\end{equation}
 \end{theorem}

 Some comments are in order. 
 
The proof of Theorems \ref{main-theorem} and \ref{masodik-fotetel} is based on a decomposition argument similar to the one carried out by Talenti \cite{Talenti} and Ashbaugh and Benguria \cite{A-B} in the Euclidean framework.
In fact, we transpose the original variational problem from generic nonpositively curved spaces to the space-form $N_\kappa^n$ by assuming the validity of the $\kappa$-Cartan-Hadamard conjecture on $(M,g)$. By a fourth order ODE it turns out that $\Gamma_\kappa(\Omega^\star)$ is the smallest positive solution to the cross-product of suitable Gaussian hypergeometric functions (resp., Bessel functions) whenever $\kappa>0$ (resp., $\kappa=0$). 
   The aforementioned decomposition argument combined with certain oscillatory and asymptotic properties of the  hypergeometric function provides the proof of Theorem \ref{main-theorem}. 
 
 The dimensionality restriction $n\in \{2,3\}$  in Theorem \ref{masodik-fotetel} (and relation (\ref{hatar-vegtelen})) is needed not only for the validity of the $\kappa$-Cartan-Hadamard conjecture but also for some peculiar properties of the Gaussian hypergeometric function; similar phenomenon has been pointed out also by Ashbaugh and Benguria \cite{A-B} in the Euclidean setting for Bessel functions. In addition, the arguments in Theorem \ref{masodik-fotetel} work only for sets with sufficiently \textit{small} measure; unlike the usual Lebesgue measure in $\mathbb R^n$ (where  the scaling $\Gamma_0(B_0(L))=L^{-4}\Gamma_0(B_0(1))$ holds for every $L>0$), the inhomogeneity of the canonical measure on hyperbolic spaces requires the aforementioned volume-restriction. 
 The intuitive feeling we get that eigenfunctions corresponding to $\Gamma_g(\Omega)$ on a large domain $\Omega\subset M$ with strictly negative curvature may have large nodal domains whose symmetric rearrangements in $\mathbb H_{-\kappa^2}^n$ produce large geodesic balls and their 'joined' fundamental tone can be definitely lower than the expected $\Gamma_\kappa(\Omega^\star)$. In fact, our arguments show that Theorem \ref{masodik-fotetel} cannot be improved even if we restrict the setting to the model space-form $\mathbb H_{-\kappa^2}^n$. 
It remains an open question whether or not (\ref{A_B_egyenlotlenseg}) remains valid for arbitrarily large domains in any dimension $n\geq 4$; we notice however that some nonoptimal estimates of $\Gamma_g(\Omega)$ are also provided for any domain in high-dimensions (see \S \ref{subsection-dimenziok}).  The asymptotic property (\ref{asymptotics}) for $\kappa>0$ follows by an elegant asymptotic connection between hypergeometric and Bessel functions,  which is crucial in the proof of  (\ref{A_B_egyenlotlenseg})  and its accuracy is shown in Table \ref{table} (see \S \ref{subsection-3-dim}) for some values of $L\ll 1$. Clearly, (\ref{asymptotics}) is trivial for $\kappa=0$  since $\Gamma_0(B_0(L))=\mathfrak  h_{\nu}^4/L^4$ for every $L>0$. 
 
A natural question arises concerning the sharp estimate of the fundamental tone  on complete $n$-dimensional  Riemannian manifolds with Ricci curvature Ric$_{(M,g)}\geq k(n-1)$ for some $k\geq 0$. Some arguments based on the spherical Laplacian show that Bessel functions (when $k=0$) and Gaussian hypergeometric functions (when $k>0$) will play again crucial roles. Since the parameter range of the aforementioned special functions in the nonnegatively curved case is different from the present setting, further technicalities appear which require a deep analysis. Accordingly, we intend to come back to this problem in a forthcoming paper.


As an application of Theorem \ref{masodik-fotetel}, we consider the elliptic problem 
\[ \   \left\{ \begin{array}{lll}
\Delta_\kappa^2 u-\mu\Delta_\kappa u+\gamma u  = |u|^{p-2 }u& {\rm in} &   B_\kappa(L), \\
u\geq 0,\ u\in W_0^{2,2}(B_\kappa(L)),
\end{array}\right. \eqno{({\mathcal P})}\]
where $B_\kappa(L)\subset \mathbb H_{-\kappa^2}^2$ is a  hyperbolic disc and the range of parameters  $\mu,\gamma,p,\kappa$ and $L$ is specified below. By using variational arguments, one can prove the following result.

\begin{theorem}\label{alkalmazas}  
Let  $\mu\geq 0$, $\gamma\in \mathbb R$,  $p>2$, $\kappa>0$ and 
$0<L<\frac{2.1492}{\kappa}.$
	The following statements hold$:$
	
	\begin{itemize}
		\item[(i)] if $\mu=0$ and problem  $({\mathcal P})$ admits a nonzero solution  then $\gamma>-\Gamma_\kappa(B_k(L));$
		\item[(ii)] if $\mu>0$ and $\gamma>-\Gamma_\kappa(B_k(L))$ then problem  $({\mathcal P})$ admits a nonzero solution.
	\end{itemize}
\end{theorem}


The paper is organized as follows. 
In Section \ref{section-preliminaries} we recall/prove those notions/results which are indispensable in our study (space-forms, $\kappa$-Cartan-Hadamard conjecture, oscillatory properties of specific Gaussian hypergeometric functions). In Section \ref{section-3} we develop an Ashbaugh-Benguria-Talenti-type decomposition from curved spaces to space-forms. In Sections  \ref{McKean-section} and \ref{section-comparion} we provide a  McKean-type spectral gap estimate  (proof of Theorem \ref{main-theorem}) and comparison principles (proof of Theorem \ref{masodik-fotetel}) for fundamental tones, respectively. In Section  \ref{section-alkalmazas} we prove Theorem \ref{alkalmazas}.

\section{Preliminaries}\label{section-preliminaries}

\subsection{Space-forms.} Let $\kappa\geq 0$ and $N^n_{\kappa}$ be the $n$-dimensional space-form  with constant sectional curvature $-\kappa^2$. When $\kappa=0$, $N_\kappa^n=\mathbb R^n$ is the usual Euclidean space, while for $\kappa>0$, $N^n_{\kappa}$ is the $n$-dimensional hyperbolic space represented by the
Poincar\'e ball model $N_\kappa^n=\mathbb H^n_{-\kappa^2}=\{x\in \mathbb R^n:|x|<1\}$
endowed with the Riemannian metric $$g_{
	\kappa}(x)=(g_{ij}(x))_{i,j={1,...,n}}=p^2_\kappa(x)\delta_{ij},$$ where
$p_\kappa(x)=\frac{2}{\kappa(1-|x|^2)}.$  $(\mathbb
H^n_{-\kappa^2},g_\kappa)$ is a Cartan-Hadamard manifold with constant
sectional curvature $-\kappa^2$. If $\nabla$ and div denote the Euclidean gradient and divergence operator in $\mathbb R^n,$  the canonical volume form,  gradient and  Laplacian operator on $N^n_{\kappa}$ are 
$${\text d}v_\kappa(x) =\left\{ \begin{array}{lll}
 {\text d}x &\mbox{if} &  \kappa=0, \\
p^n_\kappa(x) {\text d}x &\mbox{if} &  \kappa> 0,
\end{array}\right. \ \ \nabla_\kappa u=\left\{ \begin{array}{lll}
\nabla u &\mbox{if} &  \kappa=0, \\
\frac{\nabla u}{p_\kappa^2}&\mbox{if} &  \kappa> 0,
\end{array}\right.
$$
and 
$$\Delta_\kappa u=\left\{ \begin{array}{lll}
\Delta u &\mbox{if} &  \kappa=0, \\
p_\kappa^{-n}{\rm div}(p_\kappa^{n-2}\nabla u)&\mbox{if} &  \kappa> 0,
\end{array}\right.
$$
respectively. 
The distance function on $N^n_{\kappa}$ is denoted by $d_\kappa$; the distance between the origin and $x\in N^n_{\kappa}$ is
given by $$d_{\kappa}(x):=d_{\kappa}(0,x)=\left\{ \begin{array}{lll}
|x| &\mbox{if} &  \kappa=0, \\
\frac{1}{\kappa}\ln\left(\frac{1+|x|}{1-|x|}\right)&\mbox{if} &  \kappa> 0.
\end{array}\right.
$$
The volume of the geodesic ball $B_\kappa(r) =\{x\in N^n_{\kappa}:d_{\kappa}(x)<r\}$ is 
\begin{equation}\label{volume-hyperbolic}
\ds V_\kappa(r):=V_\kappa(B_\kappa(r)) =n\omega_n\int_0^r {\mathbf s}_\kappa(\rho)^{n-1}{\rm d}\rho,
\end{equation}
where 
$${\mathbf s}_\kappa(\rho)=\left\{ \begin{array}{lll}
\rho &\mbox{if} &  \kappa=0, \\
\frac{\sinh(\kappa\rho)}{\kappa} &\mbox{if} &  \kappa> 0.
\end{array}\right.$$

A simple change of variables gives the following useful transformation. 

\begin{proposition}\label{change-of-variables}
	Let $\kappa\geq 0.$ For every integrable function $g:[0,L]\to \mathbb R$  with $L\geq 0$ one has $$\int_0^{L}g(s){\rm d}s=\int_{B_\kappa({r_L})}g(V_\kappa(d_\kappa(x))){\rm d} v_\kappa(x),$$
	where $r_L\geq 0$ is the unique real number verifying $V_k({r_L})=L.$
\end{proposition}

\subsection{$\kappa$-Cartan-Hadamard conjecture} \label{CH-conjecture}

Let $(M,g)$ be an  $n$-dimensional Cartan-Hadamard manifold  with  sectional  curvature  bounded  above  by $-\kappa^2$ for some $\kappa\geq 0$.  The {\it $\kappa$-Cartan-Hadamard conjecture} on $(M,g)$ (called also as the \textit{generalized Cartan-Hadamard conjecture}) states that the $\kappa$-sharp isoperimetric inequality holds on $(M,g)$, i.e.  
 for every open bounded  
$\Omega\subset M$ one has
	\begin{equation}\label{cartan-hadamard-conjecture}
A_g(\partial \Omega) \geq  A_{\kappa}(\partial B_\kappa(r)),
\end{equation}
whenever $V_g(\Omega)=V_\kappa(r)$; moreover, equality holds in (\ref{cartan-hadamard-conjecture}) if and only if $\Omega$ is isometric to $B_\kappa(r)$. Hereafter,  $A_g$ and  $A_{\kappa}$ stand for the
area on $(M,g)$ and  $N_\kappa^n$, respectively.

The $\kappa$-Cartan-Hadamard conjecture holds for \textit{every} $\kappa\geq 0$ on space-forms with constant sectional curvature $-\kappa^2$ (of any dimension), see Dinghas \cite{Dinghas}, and on  Cartan-Hadamard manifolds with  sectional  curvature  bounded  above  by $-\kappa^2$  of dimension $2$,
see  Bol \cite{Bol}, and of 
dimension $3$, see Kleiner \cite{Kleiner}. In addition, a very recent result of 
Ghomi and Spruck \cite{GS} states that the $0$-Cartan-Hadamard conjecture holds in any dimension; in dimension 4, the validity of the $0$-Cartan-Hadamard conjecture is due to Croke \cite{Croke}. In higher dimensions and for $\kappa >0,$ the conjecture is still  open; for a detailed discussion, see  Kloeckner and Kuperberg \cite{KK}. 

\subsection{Gaussian hypergeometric function} 
For $a,b,c\in \mathbb C$ ($c\neq 0,-1,-2,...$) the  Gaussian hypergeometric function is defined by  
\begin{equation}\label{F-ertelmezes}
{\bf F}(a,b;c;z)={_2F}_1(a,b;c;z)= \sum_{k\geq 0}\frac{(a)_k(b)_k}{(c)_k}\frac{z^k}{k!}
\end{equation}
on the disc $|z|<1$ and extended by analytic continuation elsewhere, where $(a)_k=\frac{\Gamma(a+k)}{\Gamma(a)}$ denotes the Pochhammer symbol. 
The corresponding differential equation to $z\mapsto {\bf F}(a,b;c;z)$ is 
\begin{equation}\label{hyper-ODE}
z(1-z)w''(z)+(c-(a+b+1)z)w'(z)-abw(z)=0.
\end{equation}
We also recall the differentiation formula
\begin{equation}\label{F-differential}
\frac{\rm d}{{\rm d}z}{\bf F}(a,b;c;z)=\frac{ab}{c}{\bf F}(a+1,b+1;c+1;z).
\end{equation}

Let $n\geq 2$ be an integer,  $K>0$ be fixed, and consider the function $$
w^K_\pm(t)={\bf F}\left(\frac{1-\sqrt{(n-1)^2\pm 4\sqrt{K}}}{2},\frac{1+ \sqrt{(n-1)^2\pm4\sqrt{K}}}{2};\frac{n}{2};-t\right),\ \ t>0.$$
%
%
%
%

\noindent The following result will be indispensable in our study. 
\begin{proposition}\label{proposition-oscillatory}
	Let $K>0$ be fixed. The following properties hold$:$ 
	\begin{itemize}
		\item[(i)] $w^K_+(t)>0$ for every $t\geq 0;$
		\item[(ii)] if $K\leq \frac{{(n-1)^4}}{16}$, then $w^K_+(t)\geq w^K_-(t)>0$ for every $t\geq 0;$
		\item[(iii)] $w^K_-$ is oscillatory on $(0,\infty)$ if and only if $K>\frac{{(n-1)^4}}{16}$.
	\end{itemize}
\end{proposition}

{\it Proof.}  
For simplicity of notation, let  $a_\pm=\frac{1-\sqrt{(n-1)^2\pm 4\sqrt{K}}}{2}$ and $b_\pm=\frac{1+\sqrt{(n-1)^2\pm 4\sqrt{K}}}{2}$.

(i) The connection formula  (15.10.11) of Olver \textit{et al.} \cite{Digital} implies
that $$w^K_+(t)=(1+t)^{-b_+}{\bf F}\left(\frac{n}{2}-a_+,b_+;\frac{n}{2};\frac{t}{1+t}\right),\ t\geq 0.$$ Due to $\frac{n}{2}-a_+>0$, $b_+>0$ and  (\ref{F-ertelmezes}), we have that $w^K_+(t)>0$ for every $t\geq 0.$

(ii) Fix  $0<K\leq \frac{{(n-1)^4}}{16}$. First, since $\frac{n}{2}-a_->0$ and $b_->0$,  the connection formula  (15.10.11) of \cite{Digital} together with (\ref{F-ertelmezes}) implies again that 
$$w^K_-(t)=(1+t)^{-b_-}{\bf F}\left(\frac{n}{2}-a_-,b_-;\frac{n}{2};\frac{t}{1+t}\right)>0, \ t>0.$$ By virtue of  (\ref{hyper-ODE}), an elementary transformation shows that $w_\pm:=w^K_\pm$ verifies the ordinary differential equation
\begin{equation}\label{ODE-modified-2}
t(t+1)w_\pm''(t)+\left(2t+\frac{n}{2}\right)w_\pm'(t)+\frac{1-(n-1)^2\mp 4\sqrt{K}}{4}w_\pm(t)=0,\ \ t>0.
\end{equation}
It turns out that (\ref{ODE-modified-2}) is equivalent to 
\begin{equation}\label{ODE-modified}
\left(p(t)
w_\pm'(t)\right)'+q_\pm (t)w_\pm(t)=0,\ \ t>0,
\end{equation}
where $p(t)=t^\frac{n}{2}(1+t)^{2-\frac{n}{2}}$,  $q_\pm(t)=\tilde K_\pm t^{\frac{n}{2}-1}(1+t)^{1-\frac{n}{2}}$ and 
$\tilde K_\pm =\frac{1-(n-1)^2\mp 4\sqrt{K}}{4}$. For any $\tau>0$, relation (\ref{ODE-modified}) and a Sturm-type argument gives  that 
\begin{eqnarray*}
0&=&\int_0^\tau\left[w_-\left(\left(p
w_+'\right)'+q_+ w_+\right)-w_+\left(\left(p
w_-'\right)'+q_- w_-\right)\right]\\&=&\int_0^\tau (q_+-q_-) w_-w_+ + \left[p\left(w_-w_+'-w_+w_-'\right)\right]\big|_{0}^{\tau}.
\end{eqnarray*}
Since $q_+<q_-$, $p(0)=0$, and $w_\pm>0$,  we necessarily have that $w_-w_+'-w_+w_-'\geq 0$ on $(0,\infty)$. In particular, $\frac{w_+}{w_-}$ is non-decreasing on $(0,\infty)$ and since $w_+(0)=w_-(0)=1$, we have that $w_+\geq w_-$ on $[0,\infty)$. 

(iii) By (ii) we have $ w^K_-(t)>0$ for every $t>0$ whenever $0<K\leq \frac{{(n-1)^4}}{16}$,  i.e.  $w^K_-$ is not oscillatory on $(0,\infty)$ for numbers $K$ belonging to this range.  

Assume now that $K> \frac{{(n-1)^4}}{16}$. 
Since $\ds\int_{\alpha}^{\infty}\frac{1}{p(t)}{\rm d}t<\infty$ for every $\alpha>0$, one can apply the result of Sugie, Kita and Yamaoka \cite[Theorem 3.1]{SKY} (see also Hille \cite{Hille}), which states that  if 
\begin{equation}\label{p-q}
p(t)q_-(t)\left(\ds\int_{t}^{\infty}\frac{1}{p(\tau)}{\rm d}\tau\right)^2\geq \frac{1}{4}\ \ {\rm for}\ \ t\gg 1,
\end{equation}
 then the function $w^K_-$ in (\ref{ODE-modified}) is oscillatory. Due to the fact that 
$$p(t)q_-(t)\left(\ds\int_{t}^{\infty}\frac{1}{p(\tau)}{\rm d}\tau\right)^2\sim \tilde K_-\ \ \ {\rm as}\ \ t\to \infty,$$ and $\tilde K_- =\frac{1-(n-1)^2+ 4\sqrt{K}}{4}>\frac{1}{4}$, 
inequality (\ref{p-q}) trivially  holds, which concludes the proof.  
 \hfill $\square$

\begin{rem}\rm Dmitrii Karp kindly pointed out that for every $\beta\geq \frac{1}{2}$ and $t>0$, the function $x\mapsto {\bf F}(\frac{1}{2}-x,\frac{1}{2}+x;\beta;-t)$ is strictly increasing on $[0,\infty)$ from which Proposition \ref{proposition-oscillatory}/(ii)  follows; his proof is based on fine properties of the hypergeometric functions $_2F_1$ and $_3F_2$, cf.  Karp \cite{Karp}.  	
\end{rem}

%


\section{Ashbaugh-Benguria-Talenti-type decomposition: from curved spaces to space-forms}\label{section-3}

Without saying explicitly throughout this section, we put ourselves into the context of Theorem \ref{main-theorem}, i.e.   we fix an $n$-dimensional $(n\geq 2)$ Cartan Hadamard manifold $(M,g)$ with sectional curvature $\textbf{K}\leq -\kappa^2\leq 0$ $(\kappa\geq 0)$, verifying the  $\kappa$-Cartan-Hadamard conjecture (see \S \ref{CH-conjecture}).  

 Let $\Omega\subset M$ be a bounded domain. Inspired by Talenti \cite{Talenti} and Ashbaugh and Benguria \cite{A-B}, we provide in this section a decomposition argument by estimating from below the fundamental tone $\Gamma_g(\Omega)$ given in (\ref{variational-charact}) by a value coming from a two-geodesic-ball minimization problem on the space-form $N_\kappa^n.$ We first state:

\begin{proposition}\label{prop-1}
	The infimum in {\rm (\ref{variational-charact})} is achieved. 
\end{proposition}
{\it Proof.} Due to Hopf-Rinow's theorem, the set $\Omega$ is relatively compact. Consequently, the Ricci curvature  is bounded from below on $\Omega$, see e.g. Bishop and Critenden \cite[p.166]{BC}, and the injectivity radius is positive on $\Omega$, see Klingenberg \cite[Proposition 2.1.10]{Kli}. By a similar argument as in Hebey \cite[Proposition 3.3]{Hebey}, based on the Bochner-Lichnerowicz-Weitzenb\"ock formula,    the norm of the Sobolev space  $H_0^{2}(\Omega)=W_0^{2,2}(\Omega)$, i.e. 
 $u\mapsto \left(\ds\int_{\Omega} (|\nabla_g^2u|^2+|\nabla_gu|^2+u^2){\rm d}v_g\right)^{1/2},$  is equivalent to the norm given by $u\mapsto \left(\ds\int_{\Omega} \left((\Delta_g u)^2+|\nabla_gu|^2+u^2\right){\rm d}v_g\right)^{1/2}.$ Accordingly, {\rm (\ref{variational-charact})}  is well-defined. The proof of the claim, i.e.  putting minimum  in {\rm (\ref{variational-charact})},  follows by a similar variational argument as in Ashbaugh and Benguria \cite[Appendix 2]{A-B}. 
\hfill $\square$\\

We are going to use certain symmetrization arguments \`a la Schwarz; 
%
namely, if $U:\Omega\to  [0,\infty)$ is a measurable function,  we introduce its equimeasurable  rearrangement function  $U^\star:N_\kappa^n\to [0,\infty)$ such that for every $t>0$ we have
\begin{equation}\label{U-definicio}
V_\kappa(\{x\in N_\kappa^n:U^\star(x)>t \})=V_g(\{x\in\Omega:U(x)>t \}).
\end{equation}
If $S\subset M$ is a measurable set, then $S^\star$ denotes the geodesic ball in $N_\kappa^n$ with center in the origin such that $V_g(S)=V_\kappa(S^\star).$

Let $u\in W_0^{2,2}(\Omega)$ be a minimizer in (\ref{variational-charact}); since $u$ is not necessarily of constant sign, let $u_+=\max(u,0)$ and $u_-=-\min(u,0)$ be the positive and negative parts of $u$, and $$\Omega_+=\{x\in \Omega: u_+(x)>0\}\ \ {\rm  and}\ \  \Omega_-=\{x\in \Omega: u_-(x)>0\},$$ respectively. For further use, let $a,b\geq 0$ such that 
\begin{equation}\label{a-b}
V_\kappa(a)=V_g(\Omega_+) \ \ {\rm and} \ \  V_\kappa(b)=V_g(\Omega_-).
\end{equation}
In particular,  $V_\kappa(a)+V_\kappa(b)=V_g(\Omega)=V_\kappa(L)$
for some $L>0.$
We define the functions $u_+^\star,u_-^\star:N_\kappa^n\to [0,\infty)$ such that for every $t>0,$ 
\begin{equation}\label{szimmetrizacio-U}
V_\kappa(\{x\in N_\kappa^n:u_+^\star(x)>t \})=V_g(\{x\in\Omega:u_+(x)>t \})=:\alpha(t),
\end{equation}
\begin{equation}\label{szimmetrizacio-U-0}
V_\kappa(\{x\in N_\kappa^n:u_-^\star(x)>t \})=V_g(\{x\in\Omega:u_-(x)>t \})=:\beta(t).
\end{equation}
The functions $u_+^\star$ and $u_-^\star$ are well-defined and radially symmetric, verifying the property that  for some $r_t>0$ and $\rho_t>0$ one has 
\begin{equation}\label{alfa-r-t}
\{x\in N_\kappa^n:u_+^\star(x)>t \}=B_\kappa(r_t)\ \ {\rm and}\ \ \{x\in N_\kappa^n:u_-^\star(x)>t \}=B_\kappa(\rho_t),\end{equation}
with $V_\kappa(r_t)=\alpha(t)$ and $V_\kappa(\rho_t)=\beta(t)$, respectively.

For further use, we consider the sets 
$$
\Lambda_t^\star=\partial(\{x\in N_\kappa^n: u_+^\star(x)>t \}),\ \ \Lambda_t=\partial(\{x\in \Omega: u_+(x)>t \}),
$$
$$
\Pi_t^\star=\partial(\{x\in N_\kappa^n: u_-^\star(x)>t \}),\ \ \Pi_t=\partial(\{x\in \Omega: -u(x)>t \}). 
$$

\begin{proposition}\label{isoperi-prop}
	Let $u\in W_0^{2,2}(\Omega)$ be a minimizer in {\rm (\ref{variational-charact})}.  Then for a.e. $t>0$ we have 
	\begin{itemize}
		\item[(i)] $\ds{A}_g(\Lambda_t)^2\leq \alpha'(t)\int_{\{u(x)>t\}}\Delta_gu{\rm d}v_g;$
		\item[(ii)] $\ds{A}_g(\Pi_t)^2\leq \beta'(t)\int_{\{u(x)<-t\}}\Delta_gu{\rm d}v_g.$
	\end{itemize}
\end{proposition}
{\it Proof.} Statements (i) and (ii) are similar to those by Talenti \cite[Appendix, p.278]{Talenti} in the Euclidean setting; for completeness, we  reproduce the proof in the curved framework. By density reasons, it is enough to consider the case when $u$ is smooth.  For $h>0$, Cauchy's inequality implies
$$\left(\frac{1}{h}\int_{t<u(x)\leq t+h}|\nabla_g u(x)|{\rm d}v_g\right)^2\leq \frac{\alpha(t)-\alpha(t+h)}{h}\frac{1}{h}\int_{t<u(x)\leq t+h}|\nabla_g u(x)|^2{\rm d}v_g.$$
When $h\to 0$, the latter relation and the co-area formula (see Chavel \cite[p.86]{Chavel}) imply that 
$$ \ds{A}_g(\Lambda_t)^2\leq -\alpha'(t)\int_{\Lambda_t}|\nabla_g u|{\rm d}\mathcal H_{n-1},$$
where $\mathcal H_{n-1}$ is the $(n-1)$-dimensional Hausdorff measure. The divergence theorem gives that 
$$\int_{\Lambda_t}|\nabla_g u|{\rm d}\mathcal H_{n-1}=-\int_{\{x\in \Omega:u_+(x)>t\}}\Delta_gu{\rm d}v_g=-\int_{\{x\in \Omega:u(x)>t\}}\Delta_gu{\rm d}v_g,$$
which concludes the proof of (i). Similar arguments hold in the proof of  (ii). 
\hfill $\square$\\


Let $$F(s)=(\Delta_g u)^\#_-(s)-(\Delta_g u)^\#_+(V_g(\Omega)-s)\ \ {\rm and}\ \ G(s)=-F(V_g(\Omega)-s),\ \ s\in [0,V_g(\Omega)],$$
where $\cdot^\#$ stands for the notation $$H^\#(s)=H^\star(x)\ \ {\rm for}\ \  s=V_\kappa(d_\kappa(x)),\ x\in \Omega.$$


\begin{proposition} \label{f-hasonlitas}
	For every $t>0$ one has that 
	\begin{itemize}
		\item[(i)] $\ds\displaystyle\int_0^{\alpha(t)} F(s){\rm d}s\geq -\int_{\{u(x)>t \}}\Delta_g u(x){\rm d}v_g(x);$
		\item[(ii)] $\ds\displaystyle\int_0^{\beta(t)} G(s){\rm d}s\geq -\int_{\{u(x)<-t \}}\Delta_g u(x){\rm d}v_g(x).$
	\end{itemize}
%
\end{proposition}

{\it Proof.} We first recall a Hardy-Littlewood-P\'olya-type inequality, i.e.   if $U:\Omega\to [0,\infty)$ is an integrable function and $U^\star$ is defined by (\ref{U-definicio}), one has  for every measurable set $S\subseteq \Omega$ that
\begin{equation}\label{Hardy-Littlewood}
\int_SU{\rm }{\rm d}v_g\leq \int_{S^\star}U^\star{\rm }{\rm d}v_k;
\end{equation} 
moreover, if  $S=\Omega$, the equality holds in (\ref{Hardy-Littlewood}) as $U^\star$ being an equimeasurable  rearrangement of $U.$

(i)  Let $t>0$ be fixed. In order to complete the proof, we are going to show first that
 \begin{equation}\label{kell-becsles-1}
  \displaystyle\int_0^{\alpha(t)} (\Delta_g u)^\#_-(s){\rm d}s\geq \int_{\{u(x)>t \}}(\Delta_g u)_-(x){\rm d}v_g(x),
 \end{equation} 
 and 
 \begin{equation}\label{kell-becsles-2}
  \displaystyle\int_{\{u(x)>t \}}(\Delta_g u)_+(x){\rm d}v_g(x)\geq \int_0^{\alpha(t)} (\Delta_g u)^\#_+(V_g(\Omega)-s){\rm d}s.
 \end{equation}

To do this, let  $r_t>0$ be the unique real number with $V_k({r_t})=\alpha(t),$ see (\ref{alfa-r-t}). The estimate  (\ref{kell-becsles-1}) follows by Proposition \ref{change-of-variables} and inequality (\ref{Hardy-Littlewood}) as
\begin{eqnarray*}
\displaystyle\int_0^{\alpha(t)} (\Delta_g u)^\#_-(s){\rm d}s&=&\int_{B_\kappa(r_t)}(\Delta_g u)^\#_-(V_\kappa(d_\kappa(x))){\rm d} v_\kappa(x)=\int_{B_\kappa(r_t)}(\Delta_g u)^\star_-(x){\rm d} v_\kappa(x)\\&\geq & \int_{\{u(x)>t \}}(\Delta_g u)_-(x){\rm d}v_g(x),
\end{eqnarray*}
where we explored that  $\{x\in \Omega: u(x)>t \}^\star=B_\kappa(r_t)$,  following by  $V_k({r_t})=\alpha(t).$

 The proof of (\ref{kell-becsles-2})  is similar; for completeness, we provide its proof. By a change of variable and  Proposition \ref{change-of-variables} it turns out that 
\begin{eqnarray*}
	\displaystyle \int_0^{\alpha(t)} (\Delta_g u)^\#_+(V_g(\Omega)-s){\rm d}s&=& \int_0^{V_g(\Omega)} (\Delta_g u)^\#_+(s){\rm d}s- \int_0^{V_g(\Omega)-\alpha(t)} (\Delta_g u)^\#_+(s){\rm d}s \\&= & \int_{\Omega^\star}(\Delta_g u)^\#_+(V_\kappa(d_\kappa(x))){\rm d} v_\kappa(x)-\int_{B_\kappa(\tau_t)}(\Delta_g u)^\#_+(V_\kappa(d_\kappa(x))){\rm d} v_\kappa(x)\\&= & \int_{\Omega^\star}(\Delta_g u)^\star_+(x){\rm d} v_\kappa(x)-\int_{B_\kappa(\tau_t)}(\Delta_g u)^\star_+(x){\rm d} v_\kappa(x),
\end{eqnarray*}
where $\tau_t>0$ is the unique real number verifying $V_\kappa(\tau_t)=V_g(\Omega)-\alpha(t)$. Let $A_t=\{x\in \Omega:u(x)\leq t\}$; then $V_g(A_t)=V_g(\Omega)-\alpha(t)=V_\kappa(\tau_t)$. In particular, by inequality (\ref{Hardy-Littlewood}) (together with the equality for the whole domain) and the latter relations we have 
\begin{eqnarray*}
	\displaystyle \int_0^{\alpha(t)} (\Delta_g u)^\#_+(V_g(\Omega)-s){\rm d}s&= & \int_{\Omega^\star}(\Delta_g u)^\star_+(x){\rm d} v_\kappa(x)-\int_{B_\kappa(\tau_t)}(\Delta_g u)^\star_+(x){\rm d} v_\kappa(x)\\&\leq  & \int_{\Omega}(\Delta_g u)_+(x){\rm d} v_g(x)-\int_{A_t}(\Delta_g u)_+(x){\rm d} v_g(x)\\&=  &\int_{\Omega\setminus A_t}(\Delta_g u)_+(x){\rm d} v_g(x)=\int_{\{u(x)>t \}}(\Delta_g u)_+(x){\rm d}v_g(x),
\end{eqnarray*}
which concludes the proof of (\ref{kell-becsles-2}). 

By (\ref{kell-becsles-1}) and (\ref{kell-becsles-2}) one has 
\begin{eqnarray*}
\int_0^{\alpha(t)} F(s){\rm d}s&=&\int_0^{\alpha(t)} (\Delta_g u)^\#_-(s){\rm d}s-\int_0^{\alpha(t)}(\Delta_g u)^\#_+(V_g(\Omega)-s){\rm d}s\\&\geq& 
\int_{\{u(x)>t \}}(\Delta_g u)_-(x){\rm d}v_g(x)-\int_{\{u(x)>t \}}(\Delta_g u)_+(x){\rm d}v_g(x)\\
&=&-\int_{\{u(x)>t \}}\Delta_g u(x){\rm d}v_g(x),
\end{eqnarray*}
which is precisely our claim. The proof of (ii) is similar. 
\hfill $\square$\\


We consider the function   $v:\Omega^\star=B_\kappa(L)\to \mathbb R$ defined by 
\begin{equation}\label{v-definit}
v(x)=\frac{1}{n\omega_n}\int_{d_\kappa(x)}^a {\mathbf s}_\kappa(\rho)^{1-n}\left(\int_0^{V_\kappa(\rho)}F(s){\rm d}s\right){\rm d}\rho.
%
\end{equation}
 A direct computation shows that $v$ is a solution to 
the problem \begin{equation}\label{CP-problem-v}
\left\{ \begin{array}{lll}
-\Delta_\kappa v(x)=F(V_\kappa(d_\kappa(x))) &\mbox{in} &  B_\kappa(a); \\
v=0  &\mbox{on} &  \partial B_\kappa(a).
\end{array}\right.
\end{equation}
In a similar way, the function  $w:\Omega^*=B_\kappa(L)\to \mathbb R$ given by 
\begin{equation}\label{w-definit}
w(x)=\frac{1}{n\omega_n}\int_{d_\kappa(x)}^b {\mathbf s}_\kappa(\rho)^{1-n}\left(\int_0^{V_\kappa(\rho)}G(s){\rm d}s\right){\rm d}\rho
%
\end{equation}
  is a solution to 
 \begin{equation}\label{CP-problem-w}
\left\{ \begin{array}{lll}
-\Delta_\kappa w(x)=G(V_\kappa(d_\kappa(x))) &\mbox{in} &  B_\kappa(b); \\
w=0  &\mbox{on} &  \partial B_\kappa(b).
\end{array}\right.
\end{equation}
In particular, by their definitions, it  turns out that
$$v\geq 0\ \ {\rm in}\ \ B_\kappa(a)\ \ \ {\rm and}\ \ \ w\geq 0\ \ {\rm in}\ \ B_\kappa(b).$$ In fact, much precise comparisons can be said by combining the above preparatory results:

\begin{theorem} \label{talenti-result}
	Let  $v$ and $w$ from {\rm (\ref{v-definit})} and {\rm (\ref{w-definit})}, respectively. Then
	\begin{equation}\label{1-compar}
	u_+^\star\leq v\ \ {in}\ \ B_\kappa(a);
	\end{equation}
	\begin{equation}\label{2-compar}
	u_-^\star\leq w\ \ {in}\ \ B_\kappa(b),
	\end{equation}
	where $a$ and $b$ are from {\rm (\ref{a-b})}. In particular, one has
	\begin{equation}\label{u-v-w}
		\int_{\Omega}u^2 {\rm d}v_g\leq \int_{B_\kappa(a)}v^2 {\rm d}v_\kappa+\int_{B_\kappa(b)}w^2 {\rm d}v_\kappa.
	\end{equation}
	In addition, 
	\begin{equation}\label{laplace-hasonlitas}
	\int_{\Omega} (\Delta_g u)^2{\rm d}v_g=\int_{B_\kappa(a)} (\Delta_\kappa v)^2{\rm d}v_\kappa+ \int_{B_\kappa(b)} (\Delta_\kappa w)^2{\rm d}v_\kappa.
	\end{equation}
\end{theorem}

{\it Proof.} We first  prove (\ref{1-compar}).
Since $(M,g)$ verifies the $\kappa$-Cartan-Hadamard conjecture, on account of (\ref{szimmetrizacio-U}) and (\ref{szimmetrizacio-U-0}), it follows that 
\begin{equation}\label{Area-Area}
{A}_\kappa(\Lambda_t^\star)\leq {A}_g(\Lambda_t)   \ \ {\rm for\  a.e.}\ t>0,
\end{equation}
\begin{equation}\label{Area-Area-2}
{A}_\kappa(\Pi_t^\star)\leq {A}_g(\Pi_t)   \ \ {\rm for\ a.e.}\ t>0.
\end{equation}
By relation (\ref{Area-Area}) and Propositions \ref{isoperi-prop} and \ref{f-hasonlitas}, one has for a.e. $t>0$ that
$${A}_\kappa(\Lambda_t^\star)^2\leq -\alpha'(t)\int_0^{\alpha(t)} F(s){\rm d}s.$$ Due to (\ref{szimmetrizacio-U}), (\ref{alfa-r-t}) and (\ref{volume-hyperbolic}),  it follows that for a.e. $t>0$, 
$$\alpha'(t)={A}_\kappa(\Lambda_t^\star)r_t'=n\omega_n{\mathbf s}_\kappa(r_t)^{n-1}r_t'.$$
Combining the above relations, it yields 
$$n\omega_n\leq -r_t'{\mathbf s}_\kappa(r_t)^{1-n}\int_0^{V_\kappa(r_t)} F(s){\rm d}s.$$
After an integration, we obtain for every $\tau\in [0,\|u_+\|_{L^\infty(\Omega)}]$ that
$$n\omega_n\tau\leq -\int_0^\tau r_t'{\mathbf s}_\kappa(r_t)^{1-n}\int_0^{V_\kappa(r_t)} F(s){\rm d}s{\rm d}t.$$
By changing the variable $r_t=\rho$, and taking into account that $r_0=a$, it follows that
$$\tau\leq \frac{1}{n\omega_n}\int_{r_\tau}^a {\mathbf s}_\kappa(\rho)^{1-n}\left(\int_0^{V_\kappa(\rho)} F(s){\rm d}s\right){\rm d}\rho.$$
Let $x\in B_\kappa(a)$ be  arbitrarily fixed and associate to this element the unique $\tau\in [0,\|u_+\|_{L^\infty(\Omega)}]$ such that $d_\kappa(x)=r_\tau$. By the definition  of $u_+^\star$ it follows that $u_+^\star(x)=\tau,$ thus the latter inequality together with (\ref{v-definit}) implies that 
$$u_+^\star(x)\leq \frac{1}{n\omega_n}\int_{d_\kappa(x)}^a {\mathbf s}_\kappa(\rho)^{1-n}\left(\int_0^{V_\kappa(\rho)} F(s){\rm d}s\right){\rm d}\rho=v(x),$$
which is precisely the claimed relation (\ref{1-compar}).  The proof of  (\ref{2-compar}) is similar, where (\ref{Area-Area-2}) is used. 

The estimate in (\ref{u-v-w}) is immediate, since
\begin{eqnarray*}
\int_{\Omega}u^2 {\rm d}v_g&=& \int_{\Omega_+}u_+^2 {\rm d}v_g+\int_{\Omega_-}u_-^2 {\rm d}v_g=\int_{B_\kappa(a)}(u_+^\star)^2 {\rm d}v_\kappa+\int_{B_\kappa(b)}(u_-^\star)^2 {\rm d}v_\kappa\\&\leq&
\int_{B_\kappa(a)}v^2 {\rm d}v_\kappa+\int_{B_\kappa(b)}w^2 {\rm d}v_\kappa,
\end{eqnarray*}
where we apply (\ref{a-b}) together with the estimates (\ref{1-compar}) and (\ref{2-compar}), respectively. 

We now prove (\ref{laplace-hasonlitas}). On one hand, by problems (\ref{CP-problem-v}) and (\ref{CP-problem-w}), Proposition \ref{change-of-variables} and a change of variables imply  that
\begin{eqnarray*}
\int_{B_\kappa(a)} (\Delta_\kappa v)^2{\rm d}v_\kappa+ \int_{B_\kappa(b)} (\Delta_\kappa w)^2{\rm d}v_\kappa&=&\int_{B_\kappa(a)} F(V_\kappa(d_\kappa(x)))^2{\rm d}v_\kappa+ \int_{B_\kappa(b)} G(V_\kappa(d_\kappa(x)))^2{\rm d}v_\kappa\\&=& \int_{0}^{V_g(\Omega_+)} F(s)^2{\rm d}s+ \int_{0}^{V_g(\Omega_-)} G(s)^2{\rm d}s\\&=&\int_{0}^{V_g(\Omega)} F(s)^2{\rm d}s.
\end{eqnarray*}
On the other hand,
\begin{align*}
\int_{0}^{V_g(\Omega)} F(s)^2{\rm d}s
&=\int_{0}^{V_g(\Omega)} \left[(\Delta_g u)^\#_-(s)^2+(\Delta_g u)^\#_+(V_g(\Omega)-s)^2-2(\Delta_g u)^\#_-(s)(\Delta_g u)^\#_+(V_g(\Omega)-s)\right]{\rm d}s
\\
&=\int_{0}^{V_g(\Omega)} \left[(\Delta_g u)^\#_-(s)^2+(\Delta_g u)^\#_+(s)^2-2(\Delta_g u)^\#_-(s)(\Delta_g u)^\#_+(V_g(\Omega)-s)\right]{\rm d}s.
\end{align*}
The latter term in the above integral  vanishes. Indeed, fix first $0\leq s<V_g(\{x\in \Omega:\Delta_gu(x)<0 \})$ and let $t:=V_g(\Omega)-s>V_g(\{x\in \Omega:\Delta_gu(x)\geq 0 \})=V_\kappa(d_\kappa(r_0))$ for some $r_0>0$. If $t=V_\kappa(d_\kappa(x))$ for some $x\in \Omega^\star=B_\kappa(L)$ then $|x|>r_0$, i.e.  $x\notin {\rm supp}(\Delta_g u)^\star_+={\rm cl}(B_\kappa(d_\kappa(r_0)))$ thus $$(\Delta_g u)^\#_+(V_g(\Omega)-s)=(\Delta_g u)^\#_+(t)=(\Delta_g u)^\#_+(V_\kappa(d_\kappa(x)))=(\Delta_g u)^\star_+(x)=0.$$
In the case when $V_g(\{x\in \Omega:\Delta_gu(x)\leq 0 \})<s\leq V_g(\Omega)$, a similar argument yields $(\Delta_g u)^\#_-(s)=0.$ Therefore, by Proposition \ref{change-of-variables} we have 
\begin{eqnarray*}
	\int_{0}^{V_g(\Omega)} F(s)^2{\rm d}s&=&\int_{\Omega^\star} \left[(\Delta_g u)^\#_-(V_\kappa(d_\kappa(x)))^2+(\Delta_g u)^\#_+(V_\kappa(d_\kappa(x)))^2\right]{\rm d}v_\kappa(x)\\&=&\int_{\Omega^\star} \left[(\Delta_g u)^\star_-(x)^2+(\Delta_g u)^\star_+(x)^2\right]{\rm d}v_\kappa(x)
	\\&=&\int_{\Omega} \left[(\Delta_g u)_-^2(x)+(\Delta_g u)^2_+(x)\right]{\rm d}v_g(x)
	\\&=&\int_{\Omega} (\Delta_g u)^2(x){\rm d}v_g(x),
\end{eqnarray*}
which concludes the proof. 
\hfill $\square$\\

\begin{proposition}\label{Boundary-condition}
	Let  $v$ and $w$ from {\rm (\ref{v-definit})} and {\rm (\ref{w-definit})}, respectively. Then
	$$\ds\left\{ \begin{array}{lll}
	v'(a)a^{n-1}=w'(b)b^{n-1} &\mbox{if} &  \kappa=0; \\
	v'(\tanh(\frac{\kappa a}{2}))\sinh(\frac{\kappa a}{2})^{n-1}\cosh(\frac{\kappa a}{2})^{n-3}=w'(\tanh(\frac{\kappa b}{2}))\sinh(\frac{\kappa b}{2})^{n-1}\cosh(\frac{\kappa b}{2})^{n-3} &\mbox{if} &  \kappa> 0.
	\end{array}\right.$$

\end{proposition}

{\it Proof.} By the boundary condition $\frac{\partial u}{\partial \textbf{n}}=0$ on $ \partial \Omega$, the divergence theorem implies that $$\int_{\Omega}\Delta_g u{\rm d}v_g=0.$$ Therefore, the latter relation and  Proposition \ref{change-of-variables} give
\begin{eqnarray*}
0=-\int_{\Omega}\Delta_g u{\rm d}v_g&=&\int_{\Omega}(\Delta_g u)_-{\rm d}v_g-\int_{\Omega}(\Delta_g u)_+{\rm d}v_g=\int_{\Omega^\star}(\Delta_g u)_-^\star{\rm d}v_\kappa-\int_{\Omega^\star}(\Delta_g u)_+^\star{\rm d}v_\kappa\\&=&\int_{\Omega^\star}(\Delta_g u)_-^\#(V_\kappa(d_\kappa(x))){\rm d}v_\kappa(x)-\int_{\Omega^\star}(\Delta_g u)_+^\#(V_\kappa(d_\kappa(x))){\rm d}v_\kappa(x)\\&=&
\int_0^{V_g(\Omega)}(\Delta_g u)_-^\#(s){\rm d}s-\int_0^{V_g(\Omega)}(\Delta_g u)_+^\#(s){\rm d}s\\&=&\int_0^{V_g(\Omega)}(\Delta_g u)_-^\#(s){\rm d}s-\int_0^{V_g(\Omega)}(\Delta_g u)_+^\#(V_g(\Omega)-s){\rm d}s\\&=&\int_0^{V_g(\Omega)}F(s){\rm d}s.
\end{eqnarray*}
Furthermore, by Proposition \ref{change-of-variables} and problems (\ref{v-definit}) and (\ref{w-definit}) we have 
\begin{eqnarray*}
	0=\int_0^{V_g(\Omega)}F(s){\rm d}s&=&\int_0^{V_g(\Omega_+)}F(s){\rm d}s+\int_{V_g(\Omega_+)}^{V_g(\Omega)}F(s){\rm d}s= \int_0^{V_g(\Omega_+)}F(s){\rm d}s-\int_{0}^{V_g(\Omega_-)}G(s){\rm d}s\\
	&=& \int_{B_\kappa(a)}F(V_\kappa(d_\kappa(x))){\rm d}v_\kappa(x)-\int_{B_\kappa(b)}G(V_\kappa(d_\kappa(x))){\rm d}v_\kappa(x)\\
	&=& -\int_{B_\kappa(a)}\Delta_{\kappa} v(x){\rm d}v_\kappa(x)+\int_{B_\kappa(b)}\Delta_{\kappa} w(x){\rm d}v_\kappa(x).
\end{eqnarray*}
A simple computation shows that
$$\int_{B_\kappa(a)}\Delta_{\kappa} v(x){\rm d}v_\kappa(x)=n\omega_n\left\{ \begin{array}{lll}
v'(r)r^{n-1}  &\mbox{if} &  \kappa=0, \\
v'(r)(1-r^2)^{2-n}r^{n-1}  &\mbox{if} &  \kappa>0,
\end{array}\right.\ \ {\rm where}\ \ d_k(r)=a.$$
Similar facts also hold for  $w$; it remains to transform the above quantities into trigonometric terms. \hfill $\square$\\



Summing up, Theorem \ref{talenti-result} and Proposition \ref{Boundary-condition} imply that

\begin{equation}\label{minimum-ketto}
\Gamma_g(\Omega)=\min_{u\in W_0^{2,2}(\Omega)\setminus \{0\}}\frac{\displaystyle \int_{\Omega}(\Delta_g u)^2 {\rm d}v_g}{\displaystyle \int_{\Omega}u^2 {\rm d}v_g}\geq \min_{v,w}\frac{\ds\int_{B_\kappa(a)} (\Delta_\kappa v)^2{\rm d}v_\kappa+ \int_{B_\kappa(b)} (\Delta_\kappa w)^2{\rm d}v_\kappa}{\ds\int_{B_\kappa(a)}v^2 {\rm d}v_\kappa+\int_{B_\kappa(b)}w^2 {\rm d}v_\kappa},
\end{equation}
where 
\begin{equation}\label{a-b-L}
V_\kappa(a)+V_\kappa(b)=V_g(\Omega)=V_\kappa(L),
\end{equation}
and the minimum in the right hand side of (\ref{minimum-ketto}) is taken over of all pairs of radially symmetric functions with  $v\in W^{1,2}_0(B_\kappa(a))\cap W^{2,2}(B_\kappa(a))$ and $w\in W^{1,2}_0(B_\kappa(b))\cap W^{2,2}(B_\kappa(b))$, $(v,w)\neq (0,0)$, verifying the boundary condition  
\begin{eqnarray}\label{BVC}
\ds\left\{ \begin{array}{lll}
v'(a)a^{n-1}=w'(b)b^{n-1} &\mbox{if} &  \kappa=0; \\
v'(\tanh(\frac{\kappa a}{2}))\sinh(\frac{\kappa a}{2})^{n-1}\cosh(\frac{\kappa a}{2})^{n-3}=w'(\tanh(\frac{\kappa b}{2}))\sinh(\frac{\kappa b}{2})^{n-1}\cosh(\frac{\kappa b}{2})^{n-3} &\mbox{if} &  \kappa> 0.
\end{array}\right.
\end{eqnarray}
We notice that the minimum in the right hand side of (\ref{minimum-ketto}) is achieved for every pair of $(a,b)$ verifying (\ref{a-b-L}), which can be proved similarly as in Proposition \ref{prop-1}; see also Ashbaugh and Benguria \cite[Appendix 2]{A-B} for the Euclidean case.

\section{McKean-type spectral gap estimate: proof of (\ref{lambda-becsles})}\label{McKean-section}
In this section we deal with a  McKean-type lower estimate of the two-geodesic-ball minimization value
\begin{equation}\label{R-minimizer}
R_{\nu,a,b}^\kappa:=\min_{v,w}\frac{\ds\int_{B_\kappa(a)} (\Delta_\kappa v)^2{\rm d}v_\kappa+ \int_{B_\kappa(b)} (\Delta_\kappa w)^2{\rm d}v_\kappa}{\ds\int_{B_\kappa(a)}v^2 {\rm d}v_\kappa+\int_{B_\kappa(b)}w^2 {\rm d}v_\kappa},
\end{equation}
subject to the conditions (\ref{a-b-L}) and (\ref{BVC}), respectively, where  $v\in W^{1,2}_0(B_\kappa(a))\cap W^{2,2}(B_\kappa(a))$ and $w\in W^{1,2}_0(B_\kappa(b))\cap W^{2,2}(B_\kappa(b))$ are radially functions,  $(v,w)\neq (0,0)$. 

Since (\ref{lambda-becsles}) is trivial for $\kappa=0$, we concern with the case  $\kappa>0$.  Let $a,b\geq 0$  verifying the constraint (\ref{a-b-L}) and  $$\alpha:=\sinh^2\left(\frac{\kappa a}{2}\right),\ \beta:=\sinh^2\left(\frac{\kappa b}{2}\right).$$ In terms of $\alpha$ and $\beta$, relation (\ref{a-b-L}) 
can be rewritten into 
  \begin{equation}\label{alfbet}
  	\int_0^{\frac{2}{\kappa}\sinh^{-1}(\sqrt{\alpha})} \sinh(\kappa\rho)^{n-1}{\rm d}\rho + \int_0^{\frac{2}{\kappa}\sinh^{-1}(\sqrt{\beta})} \sinh(\kappa\rho)^{n-1}{\rm d}\rho =\int_0^{L} \sinh(\kappa\rho)^{n-1}{\rm d}\rho.
  	\end{equation}
 

For simplicity of notation, let $$\lambda^4:=\lambda(\nu,\kappa,\alpha,\beta)^4=R_{\nu,a,b}^\kappa>0,$$ 
\begin{equation}\label{lambdak-pm}
\Lambda_{\pm}=\Lambda_{\pm}(\lambda,\kappa,n):=\sqrt{(n-1)^2\pm 4\frac{\lambda^2}{\kappa^2}}\in \mathbb C,
\end{equation}
and consider the functions
$$\mathcal G_{\pm}(\nu,\lambda,t):={\bf F}\left(\frac{1-\Lambda_{\pm}}{2},\frac{1+\Lambda_{\pm}}{2};\frac{n}{2};-t\right),\ t\geq 0,$$
\begin{equation}\label{K-defincio}
\mathcal K_\nu(\lambda,t):= \frac{\mathcal G'_-(\nu,\lambda,t)}{\mathcal G_-(\nu,\lambda,t)}-\frac{\mathcal G'_{+}(\nu,\lambda,t)}{\mathcal G_{+}(\nu,\lambda,t)},\ t\geq 0,
\end{equation}
respectively, where  $\mathcal G'_\pm(\nu,\lambda,t)=\frac{\rm d}{{\rm d} t}\mathcal G_\pm(\nu,\lambda,t).$\\


\begin{proposition}\label{lemma-fontos-1} For every $\alpha,\beta\geq 0$ verifying {\rm (\ref{alfbet})}, $\lambda=\lambda(\nu,\kappa,\alpha,\beta)$ fulfills the equation 
	\begin{equation}\label{determinant}
	(1+\alpha)^{\nu+1}\alpha^{\nu+1}\mathcal K_\nu(\lambda,\alpha)+(1+\beta)^{\nu+1}\beta^{\nu+1}\mathcal K_\nu(\lambda,\beta)=0.
	\end{equation}
	 Moreover, 
	\begin{equation}\label{lambda-minel-nagyobb}
	\lambda=\lambda(\nu,\kappa,\alpha,\beta)>\frac{n-1}{2}\kappa.
	\end{equation}
\end{proposition}
{\it Proof.} We prove relation (\ref{determinant}) by  splitting the proof into two parts. 

\textit{Case 1}: $\alpha\beta>0$. 
 Let $(v,w)$ be the minimizer in (\ref{R-minimizer}) for $R^\kappa_{\nu,a,b}=\lambda(\nu,\kappa,\alpha,\beta)^4=\lambda^4$; by the Euler-Lagrange equations and divergence theorem one obtains 
\begin{eqnarray}\label{Euler-Lagrange}
0&=&\nonumber\int_{B_\kappa(a)}(\Delta^2_\kappa v-\lambda^4 v)\phi{\rm d}v_\kappa + \int_{B_\kappa(b)}(\Delta^2_\kappa w-\lambda^4 w)\psi{\rm d}v_\kappa \\&&+ \int_{\partial B_\kappa(a)}\Delta_\kappa v p_\kappa^{n-2}\langle \nabla\phi,  \textbf{n}\rangle {\rm d}\sigma +  \int_{\partial B_\kappa(b)}\Delta_\kappa w p_\kappa^{n-2}\langle \nabla\psi,  \textbf{n}\rangle {\rm d}\sigma,
\end{eqnarray}
where $\textbf{n}$ is the outer unit normal vector to the given surface,  ${\rm d}\sigma$ is the induced surface measure and  $\phi\in C^2(B_\kappa(a))$ and $\psi\in C^2(B_\kappa(b))$ are radially symmetric test functions verifying the conditions
\begin{equation}\label{1-bcs}
\phi\left(\tanh\left(\frac{\kappa a}{2}\right)\right)=\psi\left(\tanh\left(\frac{\kappa b}{2}\right)\right)=0,
\end{equation}
\begin{equation}\label{2-bcs}
\phi'\left(\tanh\left(\frac{\kappa a}{2}\right)\right)\sinh\left(\frac{\kappa a}{2}\right)^{n-1}\cosh\left(\frac{\kappa a}{2}\right)^{n-3}=\psi'\left(\tanh\left(\frac{\kappa b}{2}\right)\right)\sinh\left(\frac{\kappa b}{2}\right)^{n-1}\cosh\left(\frac{\kappa b}{2}\right)^{n-3}.
\end{equation} 
 Now, choosing first $\psi=0$ and  $\phi\in C^2_0(B_\kappa(a))$, then   $\psi\in C^2_0(B_\kappa(b))$ and $\phi=0$ in (\ref{Euler-Lagrange}), we obtain 
\begin{equation}\label{1-4-rendu}
\Delta^2_\kappa v=\lambda^4 v \  \ {\rm in}\ \ B_\kappa(a),
\end{equation}
and
\begin{equation}\label{2-4-rendu}
\Delta^2_\kappa w=\lambda^4 w \  \ {\rm in}\ \ B_\kappa(b),
\end{equation}
respectively. Usual regularity arguments imply that $v\in C^\infty(B_\kappa(a))$ and $w\in C^\infty(B_\kappa(b))$. By the radial symmetry of the functions $v,w,\phi,\psi$, it follows that 
$$\int_{\partial B_\kappa\left(a\right)}\Delta_\kappa v p_\kappa^{n-2}\langle \nabla\phi,  \textbf{n}\rangle {\rm d}\sigma=n\omega_n \Delta_\kappa v\left(\tanh\left(\frac{\kappa a}{2}\right)\right)\phi'\left(\tanh\left(\frac{\kappa a}{2}\right)\right)\sinh\left(\frac{\kappa a}{2}\right)^{n-1}\cosh\left(\frac{\kappa a}{2}\right)^{n-3},$$
and
$$\int_{\partial B_\kappa\left(b\right)}\Delta_\kappa w p_\kappa^{n-2}\langle \nabla\psi, \textbf{n}\rangle {\rm d}\sigma=n\omega_n \Delta_\kappa w\left(\tanh\left(\frac{\kappa b}{2}\right)\right)\phi'\left(\tanh\left(\frac{\kappa b}{2}\right)\right)\sinh\left(\frac{\kappa b}{2}\right)^{n-1}\cosh\left(\frac{\kappa b}{2}\right)^{n-3}.$$
By using (\ref{Euler-Lagrange}), (\ref{2-bcs})-(\ref{2-4-rendu}) and the latter relations, it turns out that
\begin{equation}\label{Delta-relation}
\Delta_\kappa v\left(\tanh\left(\frac{\kappa a}{2}\right)\right)+\Delta_\kappa w\left(\tanh\left(\frac{\kappa b}{2}\right)\right)=0.
\end{equation}
Since $v$ is radially symmetric, one has that
$$\Delta_\kappa v(x)=\kappa^2\left[\frac{(1-r^2)^2}{4}v''(r)+\frac{1-r^2}{4r}\left({(n-3)r^2+n-1}\right)v'(r)\right],\ \ r=|x|.$$
Therefore, the fourth order ordinary differential equation  (\ref{1-4-rendu}), having no singularity at the origin, has the solution 
\begin{equation}\label{v-solution}
v(r)=(1-r^2)^\nu\left[A \mathcal G_{+}\left(\nu,\lambda,\frac{r^2}{1-r^2}\right) +B\mathcal G_{-}\left(\nu,\lambda,\frac{r^2}{1-r^2}\right)\right],\ r\in [0,\tanh(\kappa a/2)],
\end{equation}
for some  $A,B\in \mathbb R.$ In a similar way, for some $C,D\in \mathbb R$,  the non-singular solution of (\ref{2-4-rendu}) is 
\begin{equation}\label{w-solution}
w(r)=(1-r^2)^\nu\left[C \mathcal G_{+}\left(\nu,\lambda,\frac{r^2}{1-r^2}\right) +D\mathcal G_{-}\left(\nu,\lambda,\frac{r^2}{1-r^2}\right)\right],\ r\in [0,\tanh(\kappa b/2)].
\end{equation}
 By construction, both functions $v$ and $w$ are nonnegative, and after a suitable rescaling we may assume that $v(0)=w(0)=1$.  
Since $v$ and $w$ vanish on $\partial B_\kappa(a)$ and $\partial B_\kappa(b)$, respectively, one has that 
\begin{equation}\label{BC-1}
A\mathcal G_{+}(\nu,\lambda,\alpha)+B\mathcal G_{-}(\nu,\lambda,\alpha)=0,
\end{equation}
and 
\begin{equation}\label{BC-2}
C\mathcal G_{+}(\nu,\lambda,\beta)+D\mathcal G_{-}(\nu,\lambda,\beta)=0.
\end{equation}
The boundary condition (\ref{BVC}) combined with (\ref{BC-1}) and (\ref{BC-2}) takes the form 
\begin{equation}\label{BC-3}
 \alpha^{\nu+1}(1+\alpha)[A\mathcal G'_+(\nu,\lambda,\alpha)+B\mathcal G'_-(\nu,\lambda,\alpha)]- \beta^{\nu+1}(1+\beta)[C\mathcal G'_+(\nu,\lambda,\beta)+D\mathcal G'_-(\nu,\lambda,\beta)]=0.
\end{equation}
By exploring the recurrence relation for the hypergeometric function, an elementary computation transforms  relation (\ref{Delta-relation}) into 
\begin{equation}\label{BC-4}
(1+\alpha)^{-\nu}[A\mathcal G_+(\nu,\lambda,\alpha)-B\mathcal G_-(\nu,\lambda,\alpha)]+ (1+\beta)^{-\nu}[C\mathcal G_+(\nu,\lambda,\beta)-D\mathcal G_-(\nu,\lambda,\beta)]=0. 
\end{equation}
In order to  have nontrivial functions $v$ and $w$, the determinant of the $4\times 4$ matrix arising from the linear homogeneous equations given by (\ref{BC-1})-(\ref{BC-4}) should be zero, which is equivalent to 
$$
(1+\alpha)^{\nu+1}\alpha^{\nu+1}\left(\frac{\mathcal G'_-(\nu,\lambda,\alpha)}{\mathcal G_-(\nu,\lambda,\alpha)}-\frac{\mathcal G'_{+}(\nu,\lambda,\alpha)}{\mathcal G_{+}(\nu,\lambda,\alpha)}\right)+(1+\beta)^{\nu+1}\beta^{\nu+1}\left(\frac{\mathcal G'_-(\nu,\lambda,\beta)}{\mathcal G_-(\nu,\lambda,\beta)}-\frac{\mathcal G'_+(\nu,\lambda,\beta)}{\mathcal G_+(\nu,\lambda,\beta)}\right)=0,$$
giving precisely relation (\ref{determinant}). 

\textit{Case 2}: $\alpha\beta=0$. Without loss of generality, we may assume  $\alpha=0$; then $\tilde L:=\beta=\sinh(\frac{\kappa L}{2})^2>0.$ In this case, one has that $v\equiv 0$, thus $A=B=0$, and a
simpler discussion than in Case 1 (which implies (\ref{BC-2}) and the second term in (\ref{BC-3})) yields that $\mathcal K_\nu(\lambda,\tilde L)=0.$\\


{\it Proof of {\rm (\ref{lambda-minel-nagyobb})}.} Let us assume the contrary of (\ref{lambda-minel-nagyobb}), i.e.   $\lambda=\lambda(\nu,\kappa,\alpha,\beta)\leq \frac{n-1}{2}\kappa.$ On the one hand, applying Proposition \ref{proposition-oscillatory}/(ii) with $K:=\frac{\lambda^4}{\kappa^4}\leq \frac{(n-1)^4}{16}$, one has that $\mathcal G_+(\nu,\lambda,\alpha)\geq \mathcal G_-(\nu,\lambda,\alpha)>0$ and $\mathcal G_+(\nu,\lambda,\beta)\geq \mathcal G_-(\nu,\lambda,\beta)>0,$ respectively. 

\textit{Case 1}: $\alpha\beta>0$.
Since $v(0)=w(0)=1$, one has by (\ref{v-solution}) and (\ref{w-solution}) that $A+B=C+D=1$. By (\ref{BC-1}), (\ref{BC-2}) and $\mathcal G_+(\nu,\lambda,\alpha)\geq \mathcal G_-(\nu,\lambda,\alpha)>0$ and $\mathcal G_+(\nu,\lambda,\beta)\geq \mathcal G_-(\nu,\lambda,\beta)>0,$ it turns out that $A<0<B$ and $C<0<D$. On the other hand, relation (\ref{BC-4})  together with (\ref{BC-1}) and (\ref{BC-2}) gives that $A(1+\alpha)^{-\nu}\mathcal G_+(\nu,\lambda,\alpha)+ C(1+\beta)^{-\nu}\mathcal G_+(\nu,\lambda,\beta)=0,$ thus we necessarily have $AC<0$, a contradiction, which concludes the proof of  (\ref{lambda-minel-nagyobb}). 

\textit{Case 2}: $\alpha\beta=0$. Since $\mathcal G_\pm$ are analytical functions, by  continuity reason  and relation (\ref{determinant}) we have at once (\ref{lambda-minel-nagyobb}) by the previous case.
 \hfill $\square$\\

{\it Proof of {\rm (\ref{lambda-becsles}).}}
Due to relations (\ref{minimum-ketto}) and (\ref{lambda-minel-nagyobb}), for  every $\alpha,\beta\geq 0$ verifying {\rm (\ref{alfbet})}, we have 
$\ds\Gamma_g(\Omega)\geq R_{\nu,a,b}^\kappa=\lambda(\nu,\kappa,\alpha,\beta)^4\geq \frac{(n-1)^4}{16}\kappa^4,$  which is precisely relation (\ref{lambda-becsles}).  \hfill $\square$

\begin{remark}\rm The proof of (\ref{hatar-vegtelen}), i.e.  the optimality of (\ref{lambda-becsles}) in the case $n\in \{2,3\}$, requires some specific properties of the hypergeometric function that are discussed in the next section; therefore, we  postpone its proof to \S \ref{hatar-helyyzettt}.  
\end{remark}

\vspace{0.1cm}
\section{Comparison principles for fundamental tones: proof of Theorem \ref{masodik-fotetel} and  (\ref{hatar-vegtelen})}\label{section-comparion}

In the first part of this section we establish a two-sided estimate for the first positive solution of the equation (\ref{determinant}), valid on generic $n$-dimensional Cartan-Hadamard manifolds (verifying the $\kappa$-Cartan-Hadamard conjecture). In the second part we prove the sharp comparison principle for fundamental tones in $2$- and $3$-dimensions (proof of Theorem \ref{masodik-fotetel}). In the third part we give the proof of (\ref{hatar-vegtelen})   while  in the last subsection we discuss the difficulties arising in high-dimensions. As before, let $\nu=\frac{n}{2}-1.$

\subsection{Generic scheme.}\label{subsection-elso}
The comparison $\Gamma_g(\Omega)\geq\Gamma_\kappa(\Omega^\star)$ in any dimension directly follows by 
\begin{equation}\label{R-inequality}
R_{\nu,a,b}^\kappa\geq R_{\nu,0, L}^\kappa,
\end{equation}
for every $a,b\geq 0$ verifying (\ref{a-b-L}). 
Indeed, once we have (\ref{R-inequality}),  by (\ref{minimum-ketto}) and (\ref{R-minimizer}) 
it follows that 
\begin{equation}\label{sharp-compar}
\Gamma_g(\Omega)\geq R_{\nu,a,b}^\kappa\geq R_{\nu,0, L}^\kappa=\Gamma_\kappa(B_\kappa(L))=\Gamma_\kappa(\Omega^\star).
\end{equation}

%

When $\kappa=0$, inequality (\ref{R-inequality}) is verified by Ashbaugh and Benguria \cite{A-B} for $n\in \{2,3\}$; moreover, $\Gamma_0(\Omega^\star)=	\frac{\mathfrak  h^4_{\nu}}{L^4}$ where $V_g(\Omega)=\omega_nL^n$ and $\mathfrak  h_{\nu}$ is the 
first positive critical point of $\frac{J_\nu}{I_\nu}$, i.e.  the first positive zero of the cross product $J_\nu I_\nu'-I_\nu J_\nu'=J_\nu I_{\nu+1}+I_\nu J_{\nu+1}$. For $n\geq 4$, inequality (\ref{R-inequality}) fails for certain choices of $a$ and $b$. 

%
%

Let $\kappa>0$ be fixed and let $ \lambda_{\nu}(\alpha,\beta)=\lambda(\nu,\kappa,\alpha,\beta)$ be the first positive zero of  $$\lambda\mapsto (1+\alpha)^{\nu+1}\alpha^{\nu+1}\mathcal K_\nu(\lambda,\alpha)+(1+\beta)^{\nu+1}\beta^{\nu+1}\mathcal K_\nu(\lambda,\beta)=:\mathcal F_\nu(\lambda,\alpha,\beta),$$  see Proposition \ref{lemma-fontos-1}, where $\alpha,\beta\geq 0$ verify {\rm (\ref{alfbet})}, and $\tilde L=\sinh(\frac{\kappa L}{2})^2>0.$ In order to prove (\ref{R-inequality}), it suffices to show that 
\begin{equation}\label{R-inequality-bis}
\lambda_\nu(\alpha,\beta)\geq \lambda_\nu(0,\tilde L).
\end{equation}

Due to (\ref{lambda-minel-nagyobb}) and Proposition \ref{proposition-oscillatory}/(iii), for every $\lambda>\frac{n-1}{2}\kappa$ the function $t\mapsto \mathcal G_-(\nu,\lambda,t)$ has infinitely many zeros;  let  $\mathfrak  g_{\nu,k}(t)$  be the $k^{\rm th}$ zero of the functions $\mathcal G_-(\nu,\cdot,t)$ and  respectively.  Thus,  $t\mapsto \mathcal K_{\nu}(\lambda,t)$ has infinitely many simple poles. 
  
\begin{figure}[H]
 	\centering
 	\includegraphics[scale=0.6]{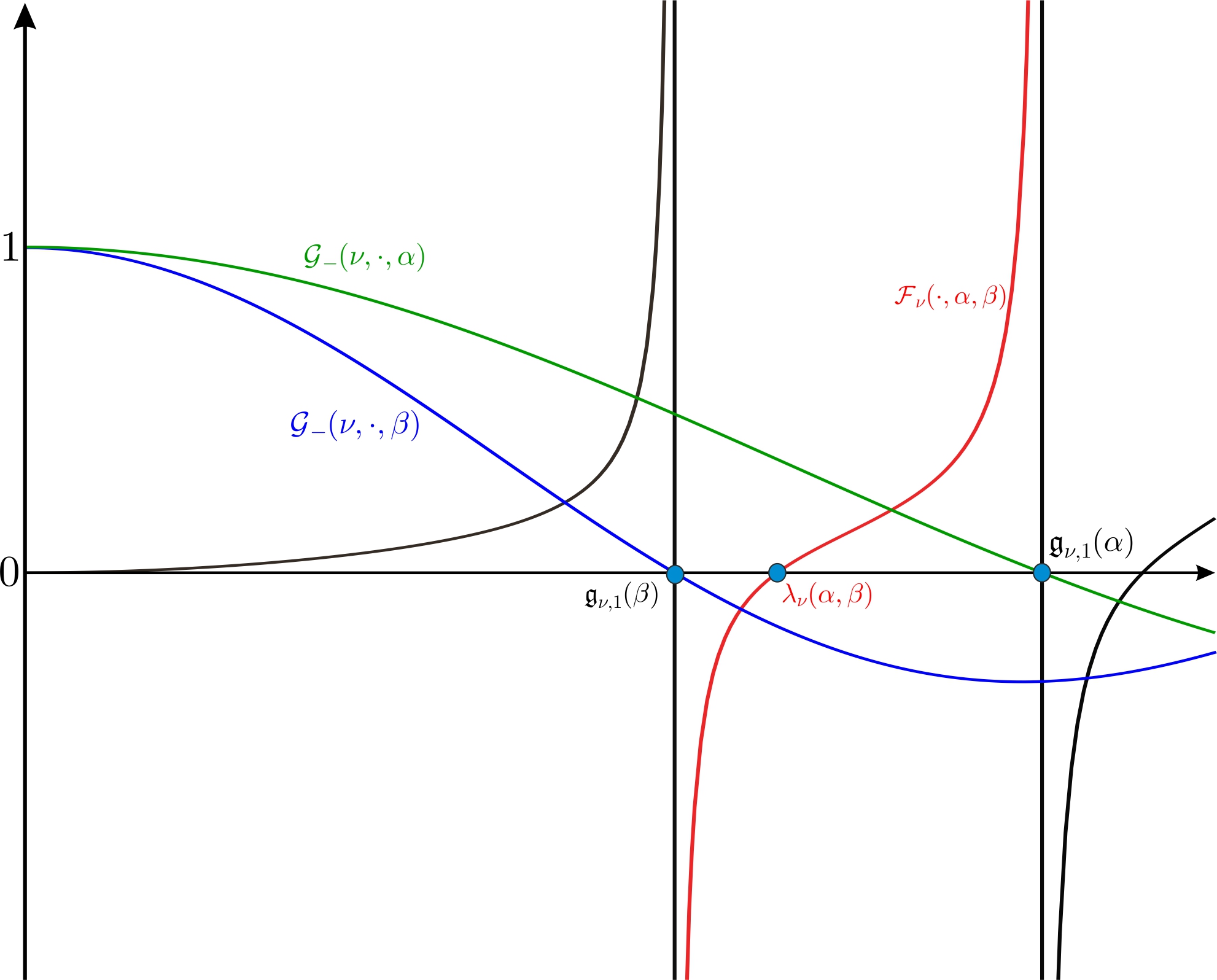}
 	\caption{The first positive zero $ \lambda_{\nu}(\alpha,\beta)$ of $\mathcal F_\nu(\cdot,\alpha,\beta)$ is between the poles $\mathfrak g_{\nu,1}(\beta)$ and $\mathfrak g_{\nu,1}(\alpha)$ of $\mathcal F_\nu(\cdot,\alpha,\beta)$; in particular, when  $\alpha$ and $\beta$ approach to $\tilde L_0=\sinh(\frac{\kappa L_0}{2})^2$ (where $2V_\kappa(L_0)=V_g(\Omega)$) it follows the limiting relation $\lambda_\nu(\tilde L_0,\tilde L_0)=\mathfrak  g_{\nu,1}(\tilde L_0).$}\label{1-abra}
 \end{figure}

Let $ L_0>0$ be fixed such that 
$2V_\kappa(L_0)=V_g(\Omega)=V_\kappa(L)$, corresponding to the case $a=b=L_0$ in 
(\ref{a-b-L}), and let $\tilde L_0=\sinh(\frac{\kappa L_0}{2})^2>0.$ 
Postponing the fact that  $\lambda\mapsto \mathcal K_\nu(\lambda,t)$ is decreasing on $(0,\infty)$ between  any two consecutive zeros of $\mathcal G_{-}(\nu,\cdot,t)$ (see Step 1 below for $\nu\in \{0,1/2\}$), and 
$\lim_{\lambda\to 0}\mathcal K_\nu(\lambda,t)=0$ for every $t>0$, the same properties hold for $\mathcal F_\nu(\cdot,\alpha,\beta)$ for any choice of $\alpha,\beta\geq 0$ verifying {\rm (\ref{alfbet})}. Accordingly,  the first positive zero $ \lambda_{\nu}(\alpha,\beta)$ of $\mathcal F_\nu(\cdot,\alpha,\beta)$ will be situated between the poles of  $\mathcal F_\nu(\cdot,\alpha,\beta)$; namely, if we assume without loss of generality that $\alpha\leq \beta,$ then 
\begin{equation}\label{g-lambda-0}
\mathfrak g_{\nu,1}(\beta)\leq \lambda_{\nu}(\alpha,\beta)\leq \min\{\mathfrak g_{\nu,1}(\alpha),\mathfrak g_{\nu,2}(\beta)\},
\end{equation}
 with the convention $\mathfrak g_{\nu,1}(0)=+\infty$. 
In the limiting case when $a$ and $b$ approach $L_0$ (thus, $\alpha$ and $\beta$ approach $\tilde L_0$), the latter relation implies that 
 $$\lambda_\nu(\tilde L_0,\tilde L_0)=\mathfrak  g_{\nu,1}(\tilde L_0),$$ 
 see  Figure \ref{1-abra}. Therefore, a necessary condition for the validity of (\ref{R-inequality-bis}) is to have 
 \begin{equation}\label{g-lambda}
 \mathfrak  g_{\nu,1}(\tilde L_0)\geq \lambda_\nu(0,\tilde L).
 \end{equation}

\begin{remark}\rm 
	Inequality (\ref{g-lambda}) \textit{fails} for every choice of $L>0$ and $\kappa\geq 0$ whenever $n\geq 4$ (thus $\nu\in \{1,3/2,2,...\}$). However, (\ref{g-lambda}) turns to be sufficient for the validity of (\ref{R-inequality-bis}) when 
	\begin{itemize}
		\item 	either $\kappa= 0$ and $n\in \{2,3\}$, corresponding to Ashbaugh and Benguria \cite{A-B};
		\item or $\kappa>0$, $n\in \{2,3\}$  and $L>0$ is \textit{sufficiently small}, see \S \ref{subsection-3-dim}. 
	\end{itemize}

\end{remark}

\subsection{The $2$- and $3$-dimensional cases: proof of Theorem \ref{masodik-fotetel}.}\label{subsection-3-dim}


In the case $\kappa=0$, relation (\ref{g-lambda}) reduces to $2^\frac{1}{n}j_{\nu,1}\geq \mathfrak  h_{\nu}$, since $\mathfrak g_{\nu,1}(\tilde L_0)=\tilde L_0^{-1}j_{\nu,1}$,  $\lambda_{\nu}(0,\tilde L_0)=\tilde L_0^{-1}\mathfrak  h_{\nu}$, and $L_0=\tilde L_0=2^{-\frac{1}{n}} L=2^{-\frac{1}{n}}\tilde L.$ Clearly, inequality $2^\frac{1}{n}j_{\nu,1}\geq \mathfrak  h_{\nu}$  holds only when $n\in \{2,3\}$, and (\ref{A_B_egyenlotlenseg}) immediately follows by (\ref{minimum-ketto}), (\ref{R-minimizer})  and the proof of Ashbaugh and Benguria \cite{A-B}, as we described in \S \ref{subsection-elso}. In addition, (\ref{asymptotics}) trivially holds since $\Gamma_0(B_0(L))=\frac{\mathfrak  h_{\nu}^4}{L^4}$ for every $L>0.$

 In the sequel, we assume that $\kappa>0$ and $n\in \{2,3\}$ (thus $\nu\in \{0,1/2\}$); the proof is divided into three steps.

{\underline{Step 1}}: {\textit{Monotonicity of $\mathcal K_{\nu}(\cdot,t)$ for $\nu\in \{0,1/2\}$.} We start with the case $n=3$ ($\nu=1/2$); the key observation is that for every $\Lambda, t>0$, one has 
$${\bf F}\left(\frac{1-i\Lambda }{2},\frac{1+i\Lambda }{2};\frac{3}{2};-t\right)=\frac{\sin(\Lambda \ln(\sqrt{t}+\sqrt{1+t}))}{\Lambda \sqrt{t}}$$
and 
$${\bf F}\left(\frac{1-\Lambda }{2},\frac{1+\Lambda }{2};\frac{3}{2};-t\right)=\frac{\sinh(\Lambda \ln(\sqrt{t}+\sqrt{1+t}))}{\Lambda \sqrt{t}},$$
both reduction formulas following by relation (15.4.15) of Olver \textit{et al.} \cite{Digital}. 
Taking advantage of the latter reduction forms, one has that 
$$\mathcal G_-(1/2,\lambda,t)=\left\{ \begin{array}{lll}
\frac{\sin(\tilde \Lambda_- \ln(\sqrt{t}+\sqrt{1+t}))}{\tilde \Lambda_- \sqrt{t}} &\mbox{if} &  \lambda>\kappa; \\
\frac{\ln(\sqrt{t}+\sqrt{1+t})}{\sqrt{t}} &\mbox{if} &  \lambda=\kappa; \\
\frac{\sinh(\Lambda_- \ln(\sqrt{t}+\sqrt{1+t}))}{ \Lambda_- \sqrt{t}} &\mbox{if} &  \lambda<\kappa, \\
\end{array}\right.
\ \ {\rm and}\ \ \mathcal G_+(1/2,\lambda,t)=\frac{\sinh(\tilde\Lambda_+ \ln(\sqrt{t}+\sqrt{1+t}))}{\tilde \Lambda_+ \sqrt{t}},$$
where 
\begin{equation}\label{lambdak}
\tilde \Lambda_-:=i\Lambda_-= 2\sqrt{\frac{\lambda^2}{\kappa^2}-1}\ \ {\rm and}\ \ \tilde \Lambda_+:=\Lambda_+=2\sqrt{\frac{\lambda^2}{\kappa^2}+1}.
\end{equation}
Thus, by (\ref{K-defincio})   one has for every $ t>0$ that
 \begin{equation}\label{explicit3dimenzioban}
\mathcal K_{1/2}(\lambda,t)={\small \frac{1}{2\sqrt{t(1+t)}}\cdot
\left\{ \begin{array}{lll}
{\tilde \Lambda_-\cot(\tilde \Lambda_-\ln(\sqrt{t}+\sqrt{1+t}))-\tilde \Lambda_+\coth(\tilde \Lambda_+\ln(\sqrt{t}+\sqrt{1+t}))} &\mbox{if} &  \lambda>\kappa; \\
\frac{1}{\ln(\sqrt{t}+\sqrt{1+t})}-2\sqrt{2}\coth(2\sqrt{2}\ln(\sqrt{t}+\sqrt{1+t})) &\mbox{if} &  \lambda=\kappa; \\
{\Lambda_-\coth( \Lambda_-\ln(\sqrt{t}+\sqrt{1+t}))-\tilde \Lambda_+\coth(\tilde \Lambda_+\ln(\sqrt{t}+\sqrt{1+t}))}  &\mbox{if} &  \lambda<\kappa. \\
\end{array}\right.}
\end{equation}
Elementary computation guarantees that $\lambda\mapsto \mathcal K_{1/2}(\lambda,t)$ is decreasing on $(0,\infty)$ between any two consecutive zeros of $\mathcal G_{-}(1/2,\cdot,t)$ for every $t>0$ fixed; the zeros of $\mathcal G_{-}(1/2,\cdot,t)$ occur only beyond the value $\kappa$ and  can be explicitly given  by  
\begin{equation}\label{g-sajatertekek}
\mathfrak  g_{1/2,k}(t)=\kappa\sqrt{1+\left(\frac{k\pi}{2\ln(\sqrt{t}+\sqrt{1+t})}\right)^2},\ \ k\in \mathbb N.
\end{equation}
In addition, since $\Lambda_-(0)=\Lambda_+(0)=2$, we also have $\lim_{\lambda\to 0}\mathcal K_{1/2}(\lambda,t)=0$ for every $ t>0$. 
In particular, relation (\ref{g-lambda-0})  is justified for $\nu=1/2$.  

When $n=2$, the differentiation formula (\ref{F-differential}) and the connection formula  (15.10.11) of Olver \textit{et al.} \cite{Digital} together with (\ref{K-defincio}) give
\begin{equation}\label{explicit2dimenzioban}
\mathcal K_{0}(\lambda,t)=-\frac{\lambda^2}{\kappa^2(1+t)}\left(\frac{{\bf F}\left(\frac{1+\Lambda_{+}}{2},\frac{3+\Lambda_{+}}{2};2;\frac{t}{1+t}\right)}{{\bf F}\left(\frac{1+\Lambda_{+}}{2},\frac{1+\Lambda_{+}}{2};1;\frac{t}{1+t}\right)}+\frac{{\bf F}\left(\frac{1+\Lambda_{-}}{2},\frac{3+\Lambda_{-}}{2};2;\frac{t}{1+t}\right)}{{\bf F}\left(\frac{1+\Lambda_{-}}{2},\frac{1+\Lambda_{-}}{2};1;\frac{t}{1+t}\right)}\right),\ \ \lambda,t>0,
\end{equation}
where $\Lambda_{\pm}=\sqrt{1\pm 4\frac{\lambda^2}{\kappa^2}}$ is from (\ref{lambdak-pm}). 
It is clear that $\lim_{\lambda\to 0}\mathcal K_{0}(\lambda,t)=0$ for every $ t>0$.
By using the definition (\ref{F-ertelmezes}) of the hypergeometric functions and the continued fraction representation, see  Cuyt \textit{et al.} \cite[Chapter 15]{Cuyt}(15.7.5) and Olver \textit{et al.} \cite{Digital}, a long computation shows that for every fixed $t>0$ the function  $\lambda\mapsto \mathcal K_{0}(\lambda,t)$ is decreasing on $(0,\infty)$ between  any two consecutive zeros of $\mathcal G_{-}(0,\cdot,t)$; see also Karp \cite{Karp}. 

{\underline{Step 2}}: {\textit{Admissible range for $L>0$ in} (\ref{g-lambda}).} 
 We are going to prove that (\ref{g-lambda}) holds for small $L>0$. We first give a crucial asymptotic estimate for $\lambda_{\nu}(0,\tilde L)$ when $L\ll 1$, i.e.  assume that  
 \begin{equation}\label{aszimptota-1}
 {\lambda_{\nu}(0,\tilde L)}\sim \kappa\sqrt{\frac{(n-1)^2}{4}+\frac{C^2}{L^2}}\ \ {\rm as}\ \  L\to 0,
 \end{equation}
 for some $C>0,$ where $\tilde L=\sinh(\frac{\kappa L}{2})^2$; our computations are valid for every $\nu\in \{0,1/2,1,...\}$.
 We observe that for every $k\in \mathbb N$ one has
 $$\lim_{L\to 0}\left(\frac{1}{2}-i\frac{C}{L}\right)_k\left(\frac{1}{2}+i\frac{C}{L}\right)_k\sinh^{2k}\left(\frac{\kappa L}{2}\right)=\left(\frac{C\kappa}{2}\right)^{2k}.$$
 Thus, by (\ref{aszimptota-1}) and uniform-convergence reasons, the latter limit implies that \\
 \\
 $\lim_{L\to 0}\mathcal G_{-}(\nu,\lambda_{\nu}(0,\tilde L),\tilde L)=$
 \begin{eqnarray*}
  &=&\lim_{L\to 0}{\bf F}\left(\frac{1-i\sqrt{4\frac{{\lambda^2_{\nu}(0,\tilde L)}}{\kappa^2}-
  		(n-1)^2 }}{2},\frac{1+i\sqrt{4\frac{{\lambda^2_{\nu}(0,\tilde L)}}{\kappa^2}-
  		(n-1)^2 }}{2};\frac{n}{2};-\sinh^{2}\left(\frac{\kappa L}{2}\right)\right)\\&=& \sum_{k\geq 0}\frac{(-1)^k}{k!\left(\frac{n}{2}\right)_k}\left(\frac{C\kappa}{2}\right)^{2k}\\&=&
  	\frac{\Gamma\left(\frac{n}{2}\right)}{\left(\frac{C\kappa}{2}\right)^{\nu}}J_\nu(C\kappa).
 \end{eqnarray*}
In a similar way, it turns out that 
 \begin{eqnarray*}
	\lim_{L\to 0}\mathcal G_{+}(\nu,\lambda_{\nu}(0,\tilde L),\tilde L)&=& \sum_{k\geq 0}\frac{1}{k!\left(\frac{n}{2}\right)_k}\left(\frac{C\kappa}{2}\right)^{2k}=
	\frac{\Gamma\left(\frac{n}{2}\right)}{\left(\frac{C\kappa}{2}\right)^{\nu}}I_\nu(C\kappa).
\end{eqnarray*}
 Moreover, the differentiation formula (\ref{F-differential}) provides 
  \begin{eqnarray*}
 	\lim_{L\to 0}L^2\mathcal G_{-}'(\nu,\lambda_{\nu}(0,\tilde L),\tilde L)= C^2
 	\frac{\Gamma\left(\frac{n}{2}+1\right)}{\left(\frac{C\kappa}{2}\right)^{\nu+1}}J_{\nu+1}(C\kappa)
 \end{eqnarray*}
 and 
 \begin{eqnarray*}
 	\lim_{L\to 0}L^2\mathcal G_{+}'(\nu,\lambda_{\nu}(0,\tilde L),\tilde L)= -C^2
 	\frac{\Gamma\left(\frac{n}{2}+1\right)}{\left(\frac{C\kappa}{2}\right)^{\nu+1}}I_{\nu+1}(C\kappa).
 \end{eqnarray*}
 Since by definition $\mathcal K_{\nu}({\lambda_{\nu}(0,\tilde L)},\tilde L)=0$, the above four limits imply that
 $$\frac{J_{\nu+1}(C\kappa)}{J_{\nu}(C\kappa)}+\frac{I_{\nu+1}(C\kappa)}{I_{\nu}(C\kappa)}=0.$$
 Accordingly, we immediately have that $C\kappa=\mathfrak  h_{\nu},$  obtaining 
 \begin{equation}\label{aszimptota-1-1}
 {\lambda_{\nu}(0,\tilde L)}\sim \sqrt{\frac{(n-1)^2}{4}\kappa^2+\frac{\mathfrak  h_{\nu}^2}{L^2}}\ \ {\rm as}\ \  L\to 0,
 \end{equation}
 which is precisely (\ref{asymptotics}). 
 
 We now provide some estimates for $\mathfrak  g_{\nu,1}(\tilde L_0)$ for $\nu\in \{0,1/2\}$ whenever $L_0\to 0$. Incidentally, it turns out that for $n=2$ ($\nu=0$), the function  
 $t\mapsto \mathcal G_{-}(0,\lambda,t):={\bf F}\left(\frac{1-\Lambda_{-}}{2},\frac{1+\Lambda_{-}}{2};1;-t\right)$ appears as the extremal in the second-order Rayleigh problem (for membranes) on the geodesic ball $B_{\kappa}(L_0)$ with the initial condition ${\bf F}\left(\frac{1-\Lambda_{-}}{2},\frac{1+\Lambda_{-}}{2};1;-\tilde L_0\right)=0$  where $\tilde L_0=\sinh(\frac{\kappa L_0}{2})^2$, see e.g. Krist\'aly \cite{Kristaly}, while the first eigenvalue  $\gamma_g(B_{\kappa}(L_0))$ corresponding to (\ref{variational-charact-1}) on $B_{\kappa}(L_0)$  is precisely $\mathfrak  g_{0,1}(\tilde L_0)$. Therefore, by Chavel \cite[p.318]{Chavel} one has that 
\begin{equation}\label{Chavel-estimate}
\mathfrak  g_{0,1}(\tilde L_0)=\gamma_g(B_{\kappa}(L_0))\sim\sqrt{\frac{1}{3}\kappa^2+\left(\frac{j_{0,1}}{ L_0}\right)^2}\ \ {\rm as}\ \  L_0\to 0.
\end{equation}
For $n=3$ (thus $\nu=1/2$), since $j_{1/2,1}=\pi$, we also have by (\ref{g-sajatertekek}) that 
\begin{equation}\label{g-konkret}
\mathfrak  g_{1/2,1}(\tilde L_0)=\sqrt{\kappa^2+\left(\frac{j_{1/2,1}}{\kappa L_0}\right)^2}\  \ \ {\rm for\ all}\ L_0>0.
\end{equation}
Recalling  $2V_\kappa(L_0)=V_\kappa(L)$, it follows that $L_0\sim 2^{-\frac{1}{n}}L$ whenever $L\ll 1.$ 
 Now, by combining these facts together with (\ref{Chavel-estimate}) and (\ref{g-konkret}), 
 it follows that 
\begin{equation}\label{1-hatarertek}
\liminf_{L\to 0}\frac{\mathfrak  g_{\nu,1}(\tilde L_0)}{\lambda_{\nu}(0,\tilde L)}=
\left\{ \begin{array}{lll}
\frac{2^{\frac{1}{2}}j_{0,1}}{\mathfrak  h_{0}}\approx  \frac{2^{\frac{1}{2}}\cdot 2.4048}{3.19622}\approx 1.064>1 &\mbox{if} &  n=2,\\
\\
\frac{2^{\frac{1}{3}}j_{1/2,1}}{\mathfrak  h_{1/2}}\approx \frac{2^{\frac{1}{3}}\pi}{3.9266}\approx  1.008>1 &\mbox{if} &  n=3,
\end{array}\right.
\end{equation}
 thus verifying (\ref{g-lambda}) for sufficiently small $L>0.$

Numerical tests show that (\ref{g-lambda}) fails for large values of $L>0$ whenever $n\in \{2,3\}$; in the sequel we provide the precise proof for $n=3$. By $2V_\kappa(L_0)=V_\kappa(L)$ we observe that $L_0\sim L-\frac{\ln 2}{2\kappa }$ whenever $L\gg 1;$ in particular, (\ref{g-konkret})  shows that $\mathfrak g_{1/2,1}(\tilde L_0)\in (\mathfrak  g_{1/2,1}(\tilde L),\mathfrak  g_{1/2,2}(\tilde L)).$  Making use of (\ref{explicit3dimenzioban}) and (\ref{g-konkret}), the latter estimate implies that 
\begin{equation}\label{2-hatarertek}
\liminf_{L\to \infty}\mathcal K_{1/2}(\mathfrak g_{1/2,1}(\tilde L_0),\tilde L)=2\left(\frac{1}{\ln(2)}-\sqrt{2}\right)\approx 0.0569 >0.
\end{equation}
 If (\ref{g-lambda}) would be true for $L\gg 1,$ relation (\ref{2-hatarertek}),  the monotonicity of $\mathcal K_{1/2}(\cdot,\tilde L)$ in the interval $(\mathfrak  g_{1/2,1}(\tilde L),\mathfrak  g_{1/2,2}(\tilde L))$, see Step 1, and the fact that ${\lambda_{1/2}(0,\tilde L)}\in (\mathfrak  g_{1/2,1}(\tilde L),\mathfrak  g_{1/2,2}(\tilde L))$, see (\ref{g-lambda-0}), imply that 
$$0<\mathcal K_{1/2}(\mathfrak g_{1/2,1}(\tilde L_0),\tilde L)\leq \mathcal K_{1/2}({\lambda_{1/2}(0,\tilde L)},\tilde L)=0,$$
a contradiction.  

We now provide the approximate threshold values of $L$ when such  turnouts occur for $n=2$ and $n=3$, respectively.  
Numerical approximations show that (\ref{g-lambda}) holds for $n=2$ whenever $0<L<\frac{2.1492}{\kappa}=:l_2$ and for $n=3$ whenever $0<L<\frac{0.719}{\kappa}=:l_3$, see Figure \ref{2-abra}. 
	\begin{figure}[H]
	\centering
	\includegraphics[scale=0.5]{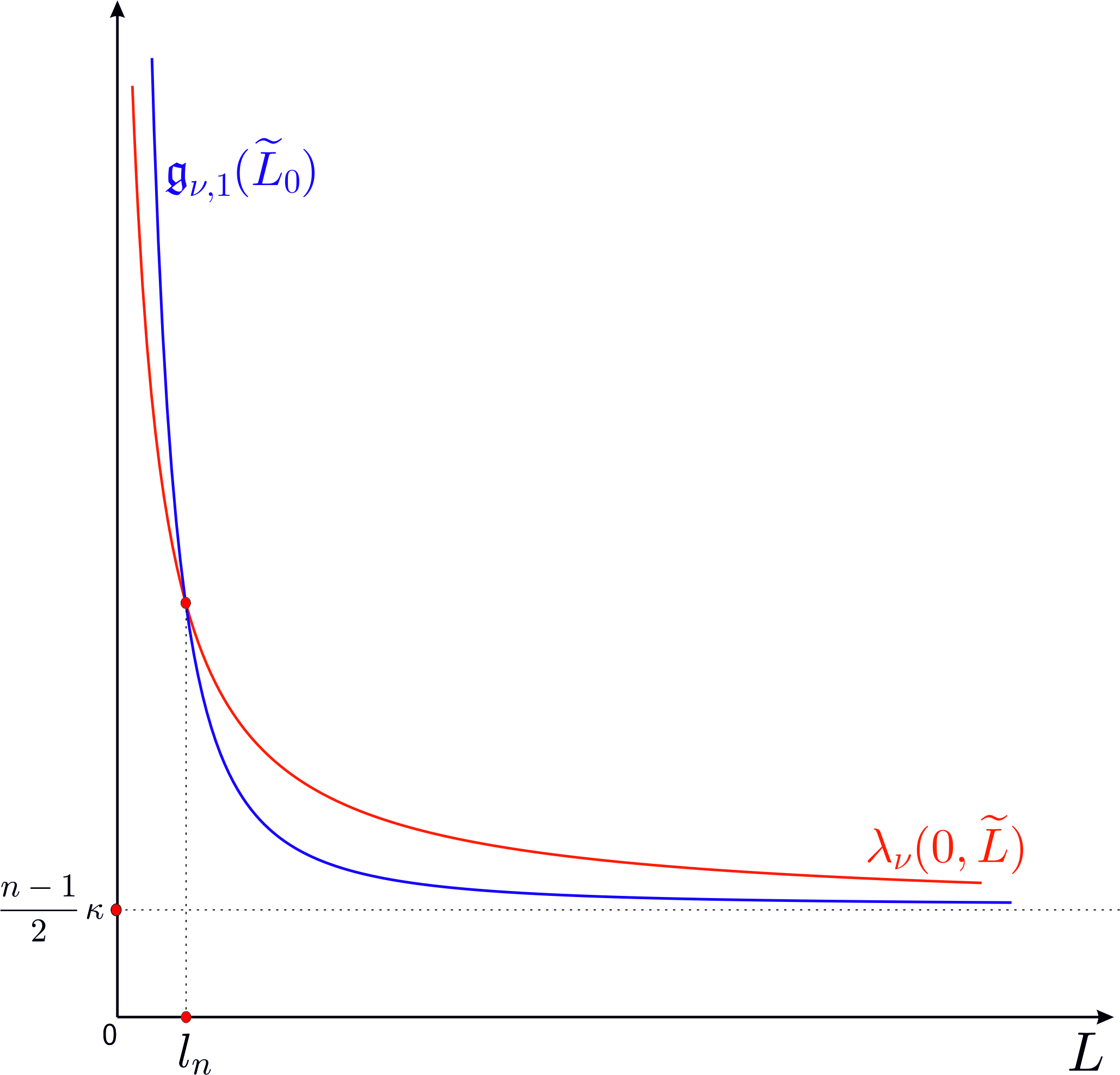}
	\caption{For $n\in \{2,3\}$ the admissible range is $0<L<l_n$ with $l_2=\frac{2.1492}{\kappa}$ and $l_3=\frac{0.719}{\kappa}$, respectively; for large values of $L$ inequality (\ref{g-lambda}) fails.}\label{2-abra}
\end{figure} 
 Due to its empirical nature, the latter values are not precise, but inequality (\ref{g-lambda}) fails for any larger values than $L=\frac{2.1493}{\kappa}$ whenever $n=2$ and $L=\frac{0.72}{\kappa}$ whenever $n=3$, respectively. Accordingly, since $V_g(\Omega)=V_\kappa(L)$, the volume of $\Omega\subset M$ cannot exceed 
\begin{equation}\label{terfogat-becsles}
V_\kappa(l_n)=n\omega_n\int_0^{l_n} {\mathbf s}_\kappa(\rho)^{n-1}{\rm d}\rho\approx 
\left\{ \begin{array}{lll}
2\pi \frac{3.34728}{\kappa^2}\approx \frac{21.031}{\kappa^2} &\mbox{if} &  n=2,\\
\\
4\pi \frac{0.137}{\kappa^3}\approx \frac{1.721}{\kappa^3} &\mbox{if} &  n=3,
\end{array}\right.
\end{equation}
which appear in the statement of the theorem.

{\underline{Step 3}}: {\textit{Concluding the proof of {\rm (\ref{A_B_egyenlotlenseg})}.} 
	 Without mentioning explicitly, we assume in the sequel that  $\alpha,\beta\geq 0$ verify {\rm (\ref{alfbet})} and  $\alpha\leq \beta$. Furthermore, without loss of generality, we may consider the case when strict inequality occurs in (\ref{g-lambda}). Since $\lambda_{\nu}(\tilde L_0,\tilde L_0)=\mathfrak  g_{{\nu},1}(\tilde L_0)>\lambda_{\nu}(0,\tilde L),$ by continuity reasons in (\ref{determinant}), it turns out that $\lambda_{{\nu}}(\alpha,\beta)>\lambda_{{\nu}}(0,\tilde L)$ for $\alpha>0$ sufficiently close to $\tilde L_0$. More precisely, the full range of $\alpha$ with this property is $[\alpha_0,\tilde L_0]$ where $\alpha_0,\beta_0$ verify {\rm (\ref{alfbet})} and $\beta_0$ is the first positive value such that $\lambda_{{\nu}}(0,\tilde L)=\mathfrak g_{{\nu},1}(\beta_0)$, i.e.  the first positive zero of $\mathcal G_-({\nu},\lambda_{{\nu}}(0,\tilde L),\cdot)$, being a pole of $\mathcal F_{{\nu}}(\lambda_{{\nu}}(0,\tilde L),\cdot,\cdot)$, see Figure \ref{3-abra}. 
	  \begin{figure}[H]
	  	\centering
	  	\includegraphics[scale=0.58]{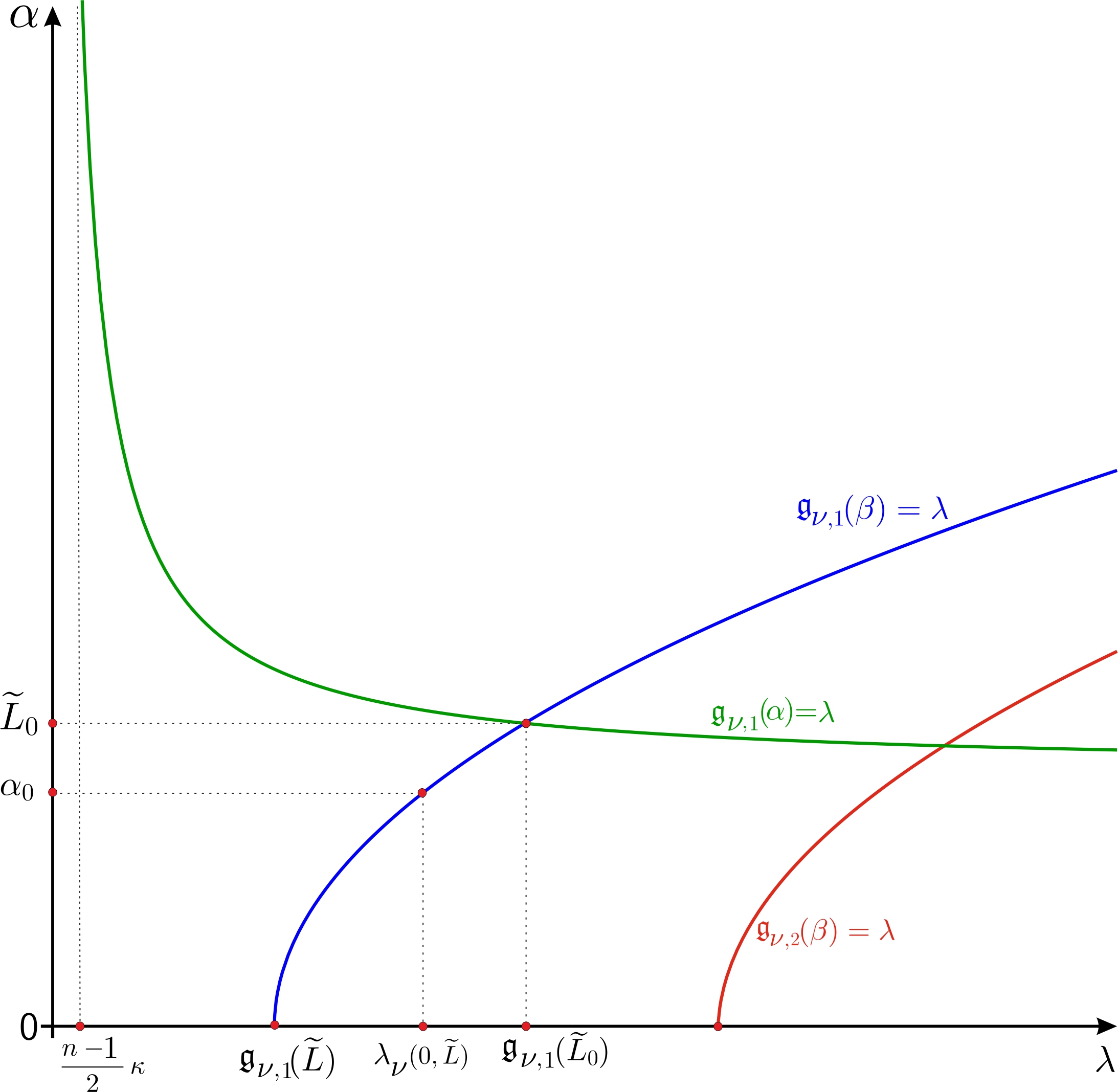}
	  	\caption{Continuity reason (when $\alpha\in [\alpha_0,\tilde L_0]$) and monotonicity argument for $\mathcal F_{\nu}$ (when $\alpha\in (0,\alpha_0)$) imply that $\lambda_{{\nu}}(\alpha,\beta)>\lambda_{{\nu}}(0,\tilde L)$. 
  		  }\label{3-abra}
	  \end{figure}  
	
We claim that for every $\alpha\in (0,\alpha_0)$, one has	
\begin{equation}\label{F-utolso}
\mathcal F_{\nu}(\lambda_{\nu}(0,\tilde L),\alpha,\beta)>0.
\end{equation}
We immediately observe that $\mathcal F_{\nu}(\lambda_{\nu}(0,\tilde L),0,\tilde L)=0$ and $\ds\lim_{\alpha\to \alpha_0^-}\mathcal F_{\nu}(\lambda_{\nu}(0,\tilde L),\alpha,\beta)=+\infty.$ In order to check (\ref{F-utolso}) one can prove that $\alpha \mapsto \mathcal F_{\nu}(\lambda_{\nu}(0,\tilde L),\alpha,\beta(\alpha))$ is increasing on $(0,\alpha_0),$ where $\beta=\beta(\alpha)$ is given by {\rm (\ref{alfbet})}. We notice that $\beta'(\alpha)=-1$ (since $\alpha+\beta=\tilde L$) when $n=2$ and  $\sqrt{\alpha(1+\alpha)}+\beta'(\alpha)\sqrt{\beta(1+\beta)}=0$ when $n=3$. Therefore, since $ \mathcal F_{\nu}$ contains ratios of hypergeometric functions, a similar monotonicity argument as in Karp and Sitnik \cite{KS} 
  implies that  $$\frac{\rm d}{\rm d \alpha} \mathcal F_{\nu}(\lambda_{\nu}(0,\tilde L),\alpha,\beta(\alpha))>0, \ \ \alpha\in (0,\alpha_0).$$    
Now, if  there exists  $\alpha\in (0,\alpha_0)$  such that  $\lambda_{{\nu}}(\alpha,\beta)< \lambda_{{\nu}}(0,\tilde L)$,  the fact that  $ \mathcal F_{\nu}(\cdot,\alpha,\beta)$ is decreasing (cf. Step 1) and relation (\ref{F-utolso}) imply that $$0<\mathcal F_{\nu}(\lambda_{\nu}(0,\tilde L),\alpha,\beta)\leq\mathcal F_{\nu}(\lambda_{\nu}(\alpha,\beta),\alpha,\beta)=0,$$ a contradiction,  
which concludes the proof of (\ref{R-inequality-bis}), so (\ref{A_B_egyenlotlenseg}).

	
	If equality occurs in (\ref{A_B_egyenlotlenseg}) then we necessarily have equality in (\ref{1-compar}) (relation (\ref{2-compar}) being canceled, or vice-versa). In particular, for a.e. $t>0$ we also have equality in  (\ref{Area-Area}), which implies equality in the isoperimetric inequality. According to the equality case in the $\kappa$-Cartan-Hadamard conjecture, the sets 
$\{x\in \Omega: u_+(x)>t \}$ and $\{x\in N_\kappa^n: u_+^\star(x)>t \}$ are isometric for a.e. $t>0$; in particular, $\Omega\subset M$ is isometric to the ball $\Omega^\star=B_\kappa(L)\subset \mathbb H_{-\kappa^2}^n$. The converse is trivial. \\

We conclude this subsection by showing the accuracy of the asymptotic estimate (\ref{asymptotics}) (see also relation (\ref{aszimptota-1-1}) in Step 2) of the fundamental tone $\Gamma_\kappa(B_\kappa(L))$  for  $L\ll 1$ in 2- and 3-dimensions; by scaling reasons, we present the values $\Gamma_\kappa(B_\kappa(L))^{1/4}$.  

\renewcommand{\arraystretch}{1.4}	
	\providecommand{\tabularnewline}{\\}
\begin{table}[h!]
	\centering
\begin{tabular}{|c|>{\centering}p{3.3cm}|>{\centering}p{3.3cm}|>{\centering}p{3.3cm}|>{\centering}p{3.3cm}|}
	\hline 
	\multirow{2}{*}{$L$} & \multicolumn{2}{c|}{$n=2$ $(\nu=0)$} & \multicolumn{2}{c|}{$n=3$ $(\nu=1/2)$}\tabularnewline
	\cline{2-5} 
	& Algebraic value of $\Gamma_\kappa(B_\kappa(L))^{1/4}$ & Approximate value of $\Gamma_\kappa(B_\kappa(L))^{1/4}$ & Algebraic value of $\Gamma_\kappa(B_\kappa(L))^{1/4}$ & Approximate value of $\Gamma_\kappa(B_\kappa(L))^{1/4}$\tabularnewline
	\hline 
		$
	0.7$ & 4.5908 &4.5728 &5.6761 & 5.6978\tabularnewline
	\hline 
	$
	0.1$ &31.9657  &31.9631 &39.2755 & 39.2787\tabularnewline
	\hline 
$	0.05$ & 63.9262&63.9248 &78.5368 & 78.5383\tabularnewline
	\hline 
	$0.003$ &1065.4069 &1065.4066 & 1308.8677& 1308.8670\tabularnewline
	\hline 
\end{tabular}
\vspace*{0.2cm}
	\caption{Comparison of the algebraic and approximate values of the fundamental tone $\Gamma_\kappa(B_\kappa(L))$ for some \textit{small} values of $L>0$; the algebraic value of $\Gamma_\kappa(B_\kappa(L))$ is  $\lambda^4$ where $\lambda>0$ is  the first positive root of $\mathcal K_\nu\left(\lambda,\sinh(\frac{\kappa L}{2})^2\right)=0$, while the approximate value of $\Gamma_\kappa(B_\kappa(L))$ is given by (\ref{asymptotics}). For simplicity,   $\kappa=1$.} \label{table}
\end{table}

\vspace{-1.0cm}
\subsection{Proof of (\ref{hatar-vegtelen}) and (\ref{magasabbbakkk})}\label{hatar-helyyzettt} We distinguish two cases. 
%

\textit{Case 1}: $n=3$. Let $L>0.$ Applying (\ref{g-lambda-0}) for $\alpha=0$ and $\beta=\tilde L=\sinh(\frac{\kappa L}{2})^2$ and using (\ref{g-sajatertekek}), it turns out that 
\begin{equation}\label{concordance}
\kappa\sqrt{1+\left(\frac{\pi}{\kappa L}\right)^2}\leq \lambda_{1/2}(0,\tilde L)\leq  \kappa\sqrt{1+\left(\frac{2\pi}{\kappa L}\right)^2}.
\end{equation}
Therefore, 
$$\lim_{L\to \infty}\Gamma_\kappa(B_\kappa(L))=\lim_{L\to \infty}\lambda_{1/2}^4(0,\tilde L)=\kappa^4,$$
which proves (\ref{hatar-vegtelen}) for $n=3$. 

 \textit{Case 2}: $n=2$. Although we have no a similar relation as  (\ref{g-sajatertekek}), we can establish its approximate version for $n=2$. We recall (see Step 2 from \S \ref{subsection-3-dim}) that the zeros  of $${\bf F}\left(\frac{1-\sqrt{1-\frac{4\lambda^2}{\kappa^2}}}{2},\frac{1+\sqrt{1-\frac{4\lambda^2}{\kappa^2}}}{2};1;-\tilde L\right)=0$$ are the values $\mathfrak  g_{0,k}(\tilde L)$, $k\in \mathbb N$,   
 and the first eigenvalue  $\gamma_g(B_{\kappa}(L))$ corresponding to (\ref{variational-charact-1}) on $B_{\kappa}(L)$  is $\mathfrak  g_{0,1}^2(\tilde L)$, 
  where $\tilde L=\sinh(\frac{\kappa L}{2})^2$. Since $\mathfrak  g_{0,k}(\tilde L)>\frac{\kappa}{2}$, see (\ref{variational-charact-1}), let $\gamma_k:=\sqrt{\frac{\mathfrak  g_{0,k}(\tilde L)^2}{\kappa^2}-\frac{1}{4}}\in \mathbb R$ and recall that $${\bf F}\left(\frac{1}{2}-i\gamma_k,\frac{1}{2}+i\gamma_k;1;-\sinh\left(\frac{\kappa L}{2}\right)^2\right)={\bf P}_{-\frac{1}{2}+i\gamma_k}(\cosh(\kappa L)),$$ where ${\bf P}_{-\frac{1}{2}+i\gamma_k}$ stands for the spherical Legendre function, see Robin \cite{Robin},  Zhurina and Karmazina \cite{ZK}. 
 By an integral representation of the spherical Legendre function, it turns  out that  for large $L>0$,  
 $$\kappa L \gamma_k\sim k\pi -\arctan(v_k/u_k),\ \ k\in \mathbb N,$$ 
 where 
 $$v_k=\int_0^\infty\frac{\sin(\gamma_k y)}{\sqrt{e^y-1}}{\rm d}y\ \ {\rm and}\ \ u_k=\int_0^\infty\frac{\cos(\gamma_k y)}{\sqrt{e^y-1}}{\rm d}y,$$
 see Zhurina and Karmazina \cite[p. 24-25]{ZK}. In particular, $\gamma_k\sim \frac{k\pi}{\kappa L}$ as $L\gg 1$ for every $k\in \mathbb N$.  Combining these facts, we have for every $k\in \mathbb N$ that 
 \begin{equation}\label{g-sajatertekek-2}
 \mathfrak  g_{0,k}(\tilde L)\sim\kappa\sqrt{\frac{1}{4}+\left(\frac{k\pi}{\kappa L}\right)^2}\ \  \ {\rm for}\ \ L\gg 1.
 \end{equation}
 By using again (\ref{g-lambda-0}) for $\alpha=0$ and $\beta=\tilde L=\sinh(\frac{\kappa L}{2})^2$, relation (\ref{g-sajatertekek-2}) provides 
 $$\kappa\sqrt{\frac{1}{4}+\left(\frac{\pi}{\kappa L}\right)^2}\leq \lambda_{0}(0,\tilde L)\leq  \kappa\sqrt{\frac{1}{4}+\left(\frac{2\pi}{\kappa L}\right)^2}\ \ \ \ {\rm for}\ \ L\gg 1.$$
 Therefore, 
 $$\lim_{L\to \infty}\Gamma_\kappa(B_\kappa(L))=\lim_{L\to \infty}\lambda_{0}^4(0,\tilde L)=\frac{\kappa^4}{16},$$
 which concludes the proof of (\ref{hatar-vegtelen}) for $n=2$.

	We now prove (\ref{magasabbbakkk}). 	In particular, by (\ref{hatar-vegtelen}) we have  for $n\in \{2,3\}$ that 
		$$\lim_{L\to \infty}\Gamma_\kappa^1(B_\kappa(L))=\lim_{L\to \infty}\Gamma_\kappa(B_\kappa(L))=\frac{(n-1)^4}{16}\kappa^4.$$  Since $\{\Gamma_\kappa^l(\Omega)\}_l$ is a nondecreasing sequence which is bounded from below by $\frac{(n-1)^4}{16}\kappa^4$ (see Proposition \ref{lemma-fontos-1}), the estimate of Cheng and Yang \cite{Cheng-Yang-1}, i.e. 
		$$ \Gamma_\kappa^{l+1}(B_\kappa(L))-\frac{(n-1)^4}{16}\kappa^4\leq 25l^{12}\left(\Gamma_\kappa^{1}(B_\kappa(L))-\frac{(n-1)^4}{16}\kappa^4\right)\ \ {\rm for\ all}\ \ l\in \mathbb N,$$
		provides  the required statement (\ref{magasabbbakkk}).  
\hfill $\square$

\begin{remark}\rm 
	(a)  The precise values of $\mathfrak  g_{1/2,k}(\tilde L)$ and the approximative values of $\mathfrak  g_{0,k}(\tilde L)$ are crucial in the proof of (\ref{hatar-vegtelen}), respectively.   The involved form of $\mathcal K_{\nu}(\lambda,t)$ for $\nu\in \{1,3/2,...\}$ (i.e.  $n\geq 4$) implies several technical difficulties to perform similar asymptotic estimates as above; however, we still believe such estimates are valid in high-dimensions.

	(b)	When $n=3$, one can give an alternative proof of  (\ref{hatar-vegtelen}). To do this, note that 
	\begin{equation}\label{ckL}
	\Gamma_\kappa(B_\kappa(L))\leq  \min_{v}\frac{\ds\int_{B_\kappa(L)} (\Delta_\kappa v)^2{\rm d}v_\kappa}{\ds\int_{B_\kappa(L)}v^2 {\rm d}v_\kappa}=:c_\kappa(L),
	\end{equation}
	where $v\in W^{2,2}_0(B_\kappa(L))\setminus \{0\}$ is taken over of all radially symmetric functions.
	A variational argument similar to the one developed in \S \ref{McKean-section} shows that   $c_\kappa(L)=\lambda^4$ where $\lambda>0$ is the first positive root  of the transcendental equation
	\begin{equation}\label{transzcendes-3-dim}
	{\tilde \Lambda_-\cot\left(\tilde \Lambda_-\frac{\kappa L}{2}\right)-\tilde \Lambda_+\coth\left(\tilde \Lambda_+\frac{\kappa L}{2}\right)}=0,
	\end{equation}
	see (\ref{explicit3dimenzioban}), where  $\tilde \Lambda_-$ and $\tilde \Lambda_+$ come from (\ref{lambdak}). 
	Analogously to (\ref{aszimptota-1}), assume that 
	\begin{equation}\label{aszimptota-2}
	c_\kappa(L)\sim \kappa^4\left({1+\frac{D^2}{L^2}}\right)^2\ \ {\rm as}\ \  L\to \infty,
	\end{equation}
	for some $D> 0.$ Inserting  (\ref{aszimptota-2})  into (\ref{transzcendes-3-dim})	and letting $L\to \infty$, a simple computation yields that $\tan(\kappa D)=0$, i.e.  $\kappa D=\pi$. We remark that (\ref{aszimptota-2}) with $D=\frac{\pi}{\kappa}$ is in a perfect concordance with  (\ref{concordance}); Table \ref{table-2} shows its accuracy (for  $\kappa=1$).  
	
	 \renewcommand{\arraystretch}{1.4}	
	 \providecommand{\tabularnewline}{\\}
	 \begin{table}[H]
	 	\centering
	 	\begin{tabular}{|c|>{\centering}p{3.4cm}|>{\centering}p{3.4cm}|>{\centering}p{3.4cm}|>{\centering}p{3.4cm}|}
	 		\hline 
	 	$L$	& Algebraic value of $\Gamma_\kappa(B_\kappa(L))^{1/4}$ & Approximate value of $\Gamma_\kappa(B_\kappa(L))^{1/4}$ \tabularnewline
	 		\hline 
	 		$
	 		50$ &1+3.1908$\cdot 10^{-3}$  &1+3.0795$\cdot 10^{-3}$\tabularnewline
	 		\hline 
	 		$	100$ & 1+5.0041$\cdot 10^{-4}$&1+4.9335$\cdot 10^{-4}$ \tabularnewline
	 		\hline 
	 		$5000$ &1+1.9745$\cdot 10^{-7}$&1+1.9739$\cdot 10^{-7}$ \tabularnewline
	 		\hline 
	 		$100000$ &1+4.71$\cdot 10^{-10}$&1+4.9348$\cdot 10^{-10}$ \tabularnewline
	 		\hline 
	 	\end{tabular}
	 	\vspace*{0.2cm}
	 	\caption{Comparison of the algebraic and approximate values of the fundamental tone $\Gamma_\kappa(B_\kappa(L))$ for some \textit{large} values of $L>0$ in $3$-dimension; the algebraic value of $\Gamma_\kappa(B_\kappa(L))$ is  $\lambda^4$ where $\lambda>0$ is  the first positive root of $\mathcal K_{1/2}\left(\lambda,\sinh(\frac{\kappa L}{2})^2\right)=0$, while the approximate value of $\Gamma_\kappa(B_\kappa(L))$ is given by (\ref{aszimptota-2}). For simplicity,  $\kappa=1$.} \label{table-2}
	 \end{table}
	 \renewcommand{\arraystretch}{1}

\end{remark}

\subsection{Fundamental tones in high-dimensions: nonoptimal estimates.}\label{subsection-dimenziok}


Our argument cannot provide sharp comparison principles for fundamental tones since inequality (\ref{g-lambda}) fails for any choice of $\kappa\geq 0$ and $L>0$ in the $n$-dimensional case whenever $n\geq 4$; we notice that similar phenomenon occurs also in the Euclidean setting, see Ashbaugh and  Benguria \cite{A-B}. However, in the case $\kappa=0$ we can provide some weak comparison principles.
To this end, if $(M,g)$ is an $n$-dimensional $(n\geq 4)$ Cartan-Hadamard manifold and $\Omega\subset M$ a bounded domain with smooth boundary, a closer inspection of the proof -- based on the validity of the $0$-Cartan-Hadamard conjecture proved by Ghomi and Spruck \cite{GS} -- gives that 	
\begin{equation}\label{CAB}
\Gamma_g(\Omega)\geq  R^0_{\nu,a,b}\geq  D_n	\Gamma_0(\Omega^\star),
\end{equation} 
where   $D_n=2^\frac{4}{n}\left(\frac{j_{\nu,1}}{\mathfrak  h_{\nu}}\right)^4$ is the constant of Ashbaugh and Laugesen \cite[Theorem 4]{A-L}.   Although  $\ds\lim_{n\to \infty}D_n=1$, the estimate (\ref{CAB}) is not sharp  since  $D_n<1$ for every $n\geq 4$. 

\vspace{0.5cm}
\section{Application: proof of Theorem \ref{alkalmazas}}\label{section-alkalmazas}

{\it Proof of} (i). 
Assume that $\mu=0$ and $({\mathcal P})$ has a nonzero solution $u\in W_0^{2,2}(B_\kappa(L))\setminus\{0\}$, i.e.  
\begin{equation}\label{diff-ujra-1}
\Delta_\kappa^2 u+\gamma u  = u^{p-1 } \ \ {\rm in}\ \ B_\kappa(L).
\end{equation}
Making use of the equation (\ref{1-4-rendu}), it turns out that the function $v(x)=v(|x|)$ given by 
$$v(r)=A \mathcal G_{+}\left(0,\lambda,\frac{r^2}{1-r^2}\right) +B\mathcal G_{-}\left(0,\lambda,\frac{r^2}{1-r^2}\right),\ r\in [0,\tanh(\kappa L/2)],$$ see  (\ref{v-solution}), is a classical solution to 
\begin{equation}\label{diff-ujra-2}
\Delta_\kappa^2 v=\Gamma_\kappa(B_k(L))v\ \ {\rm in}\ \ B_\kappa(L),
\end{equation}
while a suitable choice of the parameters $A$ and $B$ guarantees that $v\in W_0^{2,2}(B_\kappa(L))$ and  $v>0$ in $B_\kappa(L)$, respectively.  Multiplication of the equations (\ref{diff-ujra-1}) and (\ref{diff-ujra-2}) by $v>0$ and $u\geq 0$, respectively, and integrations by parts  give that 
$$\int_{B_\kappa(L)}\Delta_\kappa u\Delta_\kappa v{\rm d}v_\kappa +\gamma \int_{B_\kappa(L)}u v{\rm d}v_\kappa=\int_{B_\kappa(L)}u^{p-1} v{\rm d}v_\kappa,$$
and 
$$\int_{B_\kappa(L)}\Delta_\kappa v\Delta_\kappa u{\rm d}v_\kappa =\Gamma_\kappa(B_k(L))\int_{B_\kappa(L)} vu{\rm d}v_\kappa.$$
Therefore, one has 
$$\left(\gamma+\Gamma_\kappa(B_k(L))\right)\int_{B_\kappa(L)} vu{\rm d}v_\kappa=\int_{B_\kappa(L)}u^{p-1} v{\rm d}v_\kappa>0,$$
which immediately implies that $\gamma>-\Gamma_\kappa(B_k(L)).$

{\it Proof of} (ii). 
Let us assume that $\mu>0$ and  $\gamma>-\Gamma_\kappa(B_k(L)),$ and define the positive numbers
$$c_{\mu,\gamma}=\min\left(\mu, \frac{\min(\gamma,0)+\Gamma_\kappa(B_k(L))}{1+\Gamma_\kappa(B_k(L))}\right)\ \ {\rm and}\ \  C_{\mu,\gamma}=\max\{1,\mu,|\gamma|\}.$$ If $\mathcal T_{\mu,\gamma}:W_0^{2,2}(B_\kappa(L))\to \mathbb R$ is defined by 
$$\mathcal T_{\mu,\gamma}(u)=\int_{B_\kappa(L)} \left( (\Delta_{\kappa}u)^2+\mu|\nabla_k u|^2+\gamma u^2
\right){\rm d}v_\kappa,$$
then we have for every $u\in  W_0^{2,2}(B_\kappa(L))$ that
$$c_{\mu,\gamma}\int_{B_\kappa(L)} \left( (\Delta_{\kappa}u)^2+|\nabla_k u|^2+ u^2
\right){\rm d}v_\kappa\leq \mathcal T_{\mu,\gamma}(u)\leq C_{\mu,\gamma}\int_{B_\kappa(L)} \left( (\Delta_{\kappa}u)^2+|\nabla_k u|^2+ u^2
\right){\rm d}v_\kappa,$$
where the key ingredients are  relations (\ref{terfogat-becsles})  and  (\ref{A_B_egyenlotlenseg}), respectively.  
Therefore, $u\mapsto \mathcal T_{\mu,\gamma}^\frac{1}{2}(u)$ defines a norm on $W_0^{2,2}(B_\kappa(L))$, equivalent to the usual one, see the proof of  Proposition \ref{prop-1}. 

Let $h,H:\mathbb R\to [0,\infty)$ be defined by $h(t)=t_+^{p-1}$ and $H(t)=\frac{t_+^{p}}{p},$ where $t_+=\max(t,0)$ and associate with problem $({\mathcal P})$ its energy functional $\mathcal E_{\mu,\gamma}:W_0^{2,2}(B_\kappa(L))\to \mathbb R$ defined by
$$\mathcal E_{\mu,\gamma}(u)=\frac{1}{2}\mathcal T_{\mu,\gamma}(u)-\int_{B_\kappa(L)}H(u){\rm d}v_\kappa.$$
One can prove in a standard way that $\mathcal E_{\mu,\gamma}\in C^1(W_0^{2,2}(B_\kappa(L)); \mathbb R)$ and its differential is 
$$\mathcal E_{\mu,\gamma}'(u)(w)=\frac{1}{2}\mathcal T_{\mu,\gamma}'(u)(w)-\int_{B_\kappa(L)}h(u)w{\rm d}v_\kappa \ \ {\rm for\ all} \ \  u,w\in W_0^{2,2}(B_\kappa(L)).$$ 

We prove that $\mathcal E_{\mu,\gamma}$ satisfies the Palais-Smale condition on $W_0^{2,2}(B_\kappa(L))$. To this end, let $\{u_l\}_l\subset W_0^{2,2}(B_\kappa(L))$ be a sequence verifying $\mathcal E_{\mu,\gamma}'(u_l)\to 0$ as $l\to \infty$ and $|\mathcal E_{\mu,\gamma}(u_l)|\leq C$ ($l\in \mathbb N$) for some $C>0.$ 
The latter assumptions and relation 
$$p\mathcal E_{\mu,\gamma}(u_l)- \mathcal E_{\mu,\gamma}'(u_l)(u_l)=\frac{p}{2}\mathcal T_{\mu,\gamma}(u_l)-\frac{1}{2}\mathcal T_{\mu,\gamma}'(u_l)(u_l)\equiv \left(\frac{p}{2}-1\right)\mathcal T_{\mu,\gamma}(u_l),\ \ l\in \mathbb N,$$
immediately implies that $\{u_l\}_l$ is bounded in $W_0^{2,2}(B_\kappa(L))$; thus we may extract a  subsequence of $\{u_l\}_l$ (denoted in the same way) which weakly converges to an element $u\in W_0^{2,2}(B_\kappa(L)).$ We notice that 
$$\mathcal T_{\mu,\gamma}(u_l-u)=\mathcal E_{\mu,\gamma}'(u_l)(u_l-u)-\mathcal E_{\mu,\gamma}'(u)(u_l-u)+\int_{B_\kappa(L)}(h(u_l)-h(u))(u_l-u){\rm d}v_\kappa,\ \ l\in \mathbb N.$$
Using the fact that $\mathcal E_{\mu,\gamma}'(u_l)\to 0$ as $l\to \infty$ and $\{u_l\}_l$ is bounded in $W_0^{2,2}(B_\kappa(L))$, one has that $\mathcal E_{\mu,\gamma}'(u_l)(u_l-u)\to 0$ as $l\to \infty$. Due to the fact that $\{u_l\}_l$  weakly converges  $u$, it turns out that $\mathcal E_{\mu,\gamma}'(u)(u_l-u)\to 0$ as $l\to \infty$. Moreover, since $W_0^{2,2}(B_\kappa(L))\subset W_0^{1,2}(B_\kappa(L))\subset L^p(B_\kappa(L))$, where the latter inclusion is compact ($B_\kappa(L)\subset \mathbb H_{-\kappa^2}^2$ and $p\in (2,2^*)=(2,\infty)$), it follows that $\{u_l\}_l$ strongly converges to $u$ in $L^p(B_\kappa(L))$; therefore,  H\"older's inequality implies that 
$\ds\int_{B_\kappa(L)}(h(u_l)-h(u))(u_l-u){\rm d}v_\kappa\to 0$ as $l\to \infty$. Accordingly, $$\mathcal T_{\mu,\gamma}(u_l-u)\to 0\ \ {\rm as}\ \ l\to \infty,$$ i.e.  $\{u_l\}_l$ strongly converges to $u$ in $W_0^{2,2}(B_\kappa(L)).$

We now prove that $\mathcal E_{\mu,\gamma}$ satisfies the mountain pass geometry. First, since $p>2$, it follows that 
$$\inf_{\mathcal T_{\mu,\gamma}(u)=\rho}\mathcal E_{\mu,\gamma}(u)>0=\mathcal E_{\mu,\gamma}(0)$$
for sufficiently small  $\rho>0$.  
 Furthermore,  for sufficiently large $t>0$  and for the function $v\in W_0^{2,2}(B_\kappa(L))$ from (\ref{diff-ujra-2}) we have that 
	$$\mathcal E_{\mu,\gamma}(tv)=\frac{t^2}{2}\mathcal T_{\mu,\gamma}(v)-t^p\int_{B_\kappa(L)}H(v){\rm d}v_\kappa<0.$$

The mountain pass theorem (see e.g. Rabinowitz \cite{Rabinowitz}) implies the existence of a critical point $u\in W_0^{2,2}(B_\kappa(L))$ of $\mathcal E_{\mu,\gamma}$ with positive energy level (thus $u\neq 0$), which is nothing but a weak solution to the problem  
\[ \   \left\{ \begin{array}{lll}
\Delta_\kappa^2 u-\mu\Delta_\kappa u+\gamma u  = u_+^{p-1}&  {\rm in} &   B_\kappa(L), \\
u\in W_0^{2,2}(B_\kappa(L)).
\end{array}\right. \]
Multiplying the above equation by $u_-=\min(u,0)$, an integration on $B_\kappa(L)$ gives
$\mathcal T_{\mu,\gamma}(u_-)=0$, which implies $u_-=0$. Accordingly, $u\geq 0$   is a nonzero  solution to the original problem $({\mathcal P})$, which concludes the proof. \hfill $\square$

\begin{remark}\rm
Under the same assumptions of Theorem \ref{alkalmazas}/(ii), one can guarantee the existence of a nontrivial \textit{radially symmetric} solution to problem $({\mathcal P})$. Indeed, we can prove that the energy functional	$u\mapsto \mathcal E_{\mu,\gamma}(u)$ is invariant w.r.t. the orthogonal group $O(2)$, where the action of $O(2)$ on $W_0^{2,2}(B_\kappa(L))$ is defined by $(g*u)(x)=u(g^{-1}x)$ for every $g\in O(2)$, $x\in B_\kappa(L)$ and $u\in W_0^{2,2}(B_\kappa(L))$. Arguing in a similar way as above for the energy functional  $\mathcal E_{\mu,\gamma}^{\rm rad}:W_{0,{\rm rad}}^{2,2}(B_\kappa(L))\to \mathbb R$ instead of $\mathcal E_{\mu,\gamma}$, where $$W_{0,{\rm rad}}^{2,2}(B_\kappa(L))=\left\{u\in W_0^{2,2}(B_\kappa(L)): g*u=u\ {\rm for\ all}\ g\in O(2)\right\}$$ and $\mathcal E_{\mu,\gamma}^{\rm rad}=\mathcal E_{\mu,\gamma}\big|_{W_{0,{\rm rad}}^{2,2}(B_\kappa(L))}$, we obtain a nontrivial critical point $u_r\in W_{0,{\rm rad}}^{2,2}(B_\kappa(L))$ of $\mathcal E_{\mu,\gamma}^{\rm rad}$. Due to the principle of symmetric criticality of Palais \cite{Palais}, it turns out that  $u_r$ is a critical point of the original energy functional $\mathcal E_{\mu,\gamma}$. The rest is the same as above; moreover, since  $u_r\in W_{0,{\rm rad}}^{2,2}(B_\kappa(L))$, it follows that $u_r$ is $O(2)$-invariant, i.e.  radially symmetric. 
\end{remark}

%
%
%

\noindent {\bf Acknowledgment.} The author thanks the anonymous Referee
for her/his valuable comments and suggestions.
He is also grateful to  \'Arp\'ad Baricz, Csaba Farkas, Dmitrii Karp and Istv\'an Mez\H{o} for their help in the theory of hypergeometric functions.


\begin{thebibliography}{99}

	
	\bibitem{A-B} M. Ashbaugh, R. Benguria, On Rayleigh's conjecture for the clamped plate and its generalization to three dimensions. \textit{Duke Math. J.} 78 (1995), no. 1, 1--17.
	
	\bibitem{A-L} M. Ashbaugh, R.S. Laugesen,  Fundamental tones and buckling loads of clamped plates. \textit{Ann. Scuola Norm. Sup. Pisa Cl. Sci.} (4) 23 (1996), no. 2, 383--402.
	
\bibitem{Bauer}	L. Bauer, E. Reiss,  Block five diagonal matrices and the fast numerical solution of the biharmonic equation. \textit{Math. Comp.} 26 (1972), 311--326. 


\bibitem{BC} R.L. Bishop, R.J. Crittenden, \textit{Geometry of manifolds}.
Pure and Applied Mathematics, Vol. XV Academic Press, New York-London, 1964.


\bibitem{Bol} G. Bol,  {Isoperimetrische Ungleichungen f\"ur Bereiche auf Fl\"achen}.  \textit{Jber. Deutsch. Math. Verein.} 51 (1941), 219--257. 


\bibitem{CL1} L.M. Chasman, J.J. Langford,  On clamped plates with log-convex density. Preprint, 2018. Link:  https://arxiv.org/abs/1811.06423. 

\bibitem{CL2} L.M. Chasman, J.J. Langford, The clamped plate in Gauss space. \textit{Ann. Mat. Pura Appl.} (4)195 (2016), no. 6, 1977--2005.
	
\bibitem{Chavel} I. Chavel, \textit{Eigenvalues in Riemannian geometry}. Pure and
Applied Mathematics
115. Academic Press, Inc., Orlando, FL, 1984.

\bibitem{Chen-Zheng-Lu} D. Chen, T. Zheng, M. Lu,  Eigenvalue estimates on domains in complete noncompact Riemannian manifolds. \textit{Pacific J. Math.} 255 (2012), no. 1, 41-54.

\bibitem{Cheng-Ichikawa-Mametsuka}	Q.-M. Cheng, T. Ichikawa, S. Mametsuka, Estimates for eigenvalues of a clamped plate problem on Riemannian manifolds. \textit{J. Math. Soc. Japan}, 62 (2010), 673--686.

\bibitem{Cheng-Yang-0} Q.-M. Cheng, H. Yang,   Universal inequalities for eigenvalues of a clamped plate problem on a hyperbolic space. \textit{Proc. Amer. Math. Soc.} 139 (2011), no. 2, 461--471. 
	
	\bibitem{Cheng-Yang-1} Q.-M. Cheng, H. Yang, 
	Estimates for eigenvalues on Riemannian manifolds. 
	\textit{J. Differential Equations} 247 (2009), no. 8, 2270--2281. 
	
	
	\bibitem{Cheng-Yang-2}	Q.-M. Cheng, H. Yang, Inequalities for eigenvalues of a clamped plate problem, \textit{Trans. Amer. Math. Soc.} 358 (2006), 2625--2635.
	
\bibitem{Croke} C. Croke, A sharp four-dimensional isoperimetric inequality. \textit{Comment. Math. Helv.} 59 (2)(1984), 187--192. 

\bibitem{Coffman} C.V. Coffman,  On the structure of solutions to $\Delta^2 u = \lambda u$ which satisfy  the clamped plate conditions on a right angle. \textit{SIAM. J. Math. Anal.} 13 (1982), 746--757. 

\bibitem{CDS} C.V. Coffman, R.J. Duffin, D.H. Shaffer, \textit{The fundamental mode of vibration of a clamped annular plate is not of one sign}. Constructive approaches to mathematical models (Proc. Conf. in honor of R.J. Duffin, Pittsburgh, Pa., 1978), pp. 267--277, Academic Press, New York-London-Toronto, Ont., 1979. 
	
\bibitem{Cuyt} A. Cuyt, V.B. Petersen, B. Verdonk, H. Waadeland, W. Jones, 
	\textit{Handbook of continued fractions for special functions}.
	With contributions by Franky Backeljauw and Catherine Bonan-Hamada.  Springer, New York, 2008.
	
\bibitem{Dinghas} A. Dinghas,  Einfacher Beweis der isoperimetrischen Eigenschaft der Kugel in Riemannschen R\"aumen konstanter Kr\"ummung.  \textit{Math. Nachr.} 2 (1949), 148--162.

\bibitem{Duffin} R.J. Duffin,  Nodal lines of a vibrating plate. \textit{J. Math. Physics} 31 (1953), 294--299.

\bibitem{GS} M. Ghomi, J. Spruck, Total curvature and the isoperimetric inequality in Cartan-Hadamard manifolds, preprint, 2019. Link: https://arxiv.org/abs/1908.09814

\bibitem{Grunau} H.-C. Grunau, F. Robert,  Positivity and almost positivity of biharmonic Green's functions under Dirichlet boundary conditions. \textit{Arch. Ration. Mech. Anal.} 195 (2010), no. 3, 865--898.



	
	
%
	
	\bibitem{Hebey} E. Hebey, \textit{Nonlinear analysis on manifolds: Sobolev spaces and inequalities}. Courant Lecture Notes in Mathematics, 5.
	New York University, Courant Institute of Mathematical Sciences, New
	York; American Mathematical Society, Providence, RI, 1999.
	
	\bibitem{Hille}  E.  Hille, 
	Non-oscillation 	theorems. \textit{Trans. 	Amer. 	Math. 	Soc.}  64 (1948),  
234--252.
	 
	 \bibitem{Karp} D. Karp, private communication \& manuscript in preparation. 
	 
	  \bibitem{KS} D. Karp,  S.M. Sitnik, Inequalities and monotonicity of ratios for generalized hypergeometric function. \textit{J. Approx. Theory} 161 (2009), no. 1, 337--352.
	 
	 \bibitem{Kleiner} B. Kleiner, An isoperimetric comparison theorem. \textit{Invent. Math.} 108 (1992), no. 1, 37--–47.
	 
	 \bibitem{Kli} W. Klingenberg, \textit{Riemannian geometry}. Second edition. De Gruyter Studies in Mathematics, 1. Walter de Gruyter \& Co., Berlin, 1995.
	 
	 \bibitem{KK} 	B.R. Kloeckner, G. Kuperberg, The Cartan-Hadamard conjecture and the Little Prince, preprint, 2017. Link: https://arxiv.org/pdf/1303.3115
	 
	 \bibitem{Kristaly} A. Krist\'aly, New features of the first eigenvalue on negatively curved spaces, preprint, 2018. Link: https://arxiv.org/abs/1810.06487
	
\bibitem{McKean} H.P. McKean, An upper bound to the spectrum of $\Delta$ on a manifold of negative curvature. \textit{J. Differential Geom.} 4 (1970), 359--366. 

	\bibitem{Nad-0} N.S. Nadirashvili,  \textit{New isoperimetric inequalities in mathematical physics}. Partial differential equations of elliptic type (Cortona, 1992), 197--203, Sympos. Math., XXXV, Cambridge Univ. Press, Cambridge, 1994.

	\bibitem{Nad} N.S. Nadirashvili, Rayleigh's conjecture on the principal frequency of the clamped plate. \textit{Arch. Rational Mech. Anal.} 129 (1995), no. 1, 1--10.

		\bibitem{Digital} F.W.J. Olver, D.W. Lozier, R.F. Boisvert, C.W. Clark (eds.), \textit{NIST Handbook
	of Mathematical Functions}. Cambridge University Press, Cambridge, 2010.
	
	\bibitem{Palais} R.S. Palais, The principle of symmetric criticality. \textit{Comm. Math. Phys.} 69 (1979) 19--30. 
	
	\bibitem{Rabinowitz} P. Rabinowitz, \textit{Minimax Methods in Critical Point Theory with Applications to Differential Equations}, CBMS Regional Conference Series in Mathematics Vol. 65, Amer. Math. Soc., Providence, RI, 1986.
	
\bibitem{Rayleigh}	J.W.S. Rayleigh, \textit{The Theory of Sound}. I,II. 2nd edition, revised and enlarged. McMillan-London, London, 1926. Reprint of the 1894 and 1896 edition.
	
	\bibitem{Robin} L. Robin,  \textit{Fonctions sph\'eriques de Legendre et fonctions sph\'ero\"idales}. Tome III.  Collection Technique et Scientifique du C.N.E.T. Gauthier-Villars, Paris, 1959.
	
	\bibitem{Szego} G. Szeg\H o, On membranes and plates. \textit{Proc. Nat. Acad. Sci. U.S.A.} 36, (1950), 210--216. 
	
\bibitem{SKY} J. Sugie, K. Kita, N. Yamaoka, Oscillation constant of second-order non-linear self-adjoint differential equations. \textit{Ann. Mat. Pura Appl.} (4) 181 (2002), no. 3, 309--337.

\bibitem{Talenti} G. Talenti, 
On the first eigenvalue of the clamped plate.
\textit{Ann. Mat. Pura Appl.} (4) 129 (1981), 265--280.

\bibitem{ZK}  M.I. Zhurina, L. Karmazina,  \textit{Tables and formulae for the spherical functions} $\textbf{P}^m_{-1/2+i\tau}(z)$. Translated by E.L. Albasiny Pergamon Press, Oxford-New York-Paris, 1966.

\bibitem{Wang-Xia} Q. Wang, C. Xia, Universal bounds for eigenvalues of the biharmonic operator on Riemannian manifolds. \textit{J. Funct. Anal.} 245 (2007), 334--352.

\end{thebibliography}
\end{document}